\documentclass[a4paper,12pt,reqno,twosides,onecolumn,final]{amsart}

\usepackage{amsmath,amsfonts,amssymb,amscd}

\usepackage[foot]{amsaddr} %puts the addresses on the title page

\usepackage[skip=5pt,indent=15pt]{parskip}

\usepackage{graphicx}
\usepackage{epstopdf}
\usepackage{bm}  % \bm 
\usepackage[left,pagewise]{lineno}
\usepackage{setspace}
\usepackage[utf8]{inputenc} % set input encoding (not needed with XeLaTeX)
\usepackage{lmodern}
\usepackage{subfig}

\usepackage[]{algorithm}
\usepackage[]{algorithmic}
\usepackage{amssymb}
\usepackage{mathabx}

\usepackage{mathtools}

\usepackage{color}
\usepackage{caption}

\DeclareMathOperator{\Div}{\mathrm{div}}

\let\div\relax

\DeclareMathOperator{\div}{div}

\DeclareMathOperator{\curl}{curl}

\DeclareMathOperator{\dx}{d\Omega}

\def\0{\phantom{0}}

\newtheorem{theorem}{Theorem}[section]

\theoremstyle{definition}

\theoremstyle{remark}

\numberwithin{equation}{section}

\begin{document}


\title[Stabilized velocity post-processings for Darcy flow]{Stabilized velocity post-processings for Darcy flow in heterogeneous porous media}


%    author one information
\author[M.R. Correa]{Maicon R. Correa}
\address{Universidade Estadual de Campinas (UNICAMP), Departamento de Matemática Aplicada, IMECC, Brasil}
\curraddr{}
\email{maicon@ime.unicamp.br}
\thanks{}

%    author two information
\author[A.F.D. Loula]{Abimael F. D. Loula}
\address{LNCC - Laboratório Nacional de Computação Científica, Petrópolis, RJ, Brasil}
\curraddr{}
\email{aloc@lncc.br}
\thanks{}

%    \subjclass is required.

\date{}

\dedicatory{\rm This is a preprint of the original article published in \\{\it Communications in Numerical Methods in Engineering, 23 (2007), pp 461—489}. \\ https://onlinelibrary.wiley.com/doi/10.1002/cnm.904 \\
DOI 10.1002/cnm.904}

%    Abstract is required.
\begin{abstract}
Stable and accurate finite element methods  are presented for Darcy
flow in heterogeneous porous media with an interface of
discontinuity of the hydraulic conductivity tensor. Accurate velocity fields
are computed through  global or local post-processing formulations that use
previous approximations of the hydraulic potential. Stability is provided by
combining Galerkin and Least Squares (GLS) residuals of the governing equations with
an additional stabilization on the interface that incorporates the discontinuity
on the tangential component of the velocity field in a strong sense.
Numerical analysis is outlined and numerical results are presented
to illustrate the good performance of the proposed methods.
Convergence studies for a heterogeneous and anisotropic porous medium
confirm the same orders of convergence predicted for homogeneous
problem with smooth solutions, for both global and local post-processings.
\end{abstract}

\maketitle



%-------------------------------------------------------------
%
%
\section{Introduction}\label{sec:intro}
%
%
%-------------------------------------------------------------

Transport problems associated with reservoir simulation, groundwater
contamination or remediation demands information about the
macroscopic (average) velocity of the phases flowing in the porous
formation. These flows are, in general, highly advective and
accurate approximations of the velocity fields are crucial for the
evaluation of the related transport processes.

The system of partial differential equations that models the flow of
an incompressible homogeneous fluid in a rigid saturated porous
media consists basically of the mass conservation equation plus
Darcy's law, relating the average velocity of the fluid with the
gradient of a potential field through the hydraulic conductivity
tensor. A standard way to solve this system is based on the second
order elliptic equation obtained by substituting Darcy's law in the
mass conservation equation leading to a Poisson problem in the
potential field with velocity calculated by taking the gradient of
the solution multiplied by the hydraulic conductivity. Constructing
finite element methods based on this kind of formulation is
straightforward. However, this direct approach leads to lower-order approximations for
velocity compared to potential and, additionally, the corresponding
balance equation is enforced in an extremely weak sense. Alternative
formulations have been employed to enhance the velocity
approximation like post-processing techniques and  mixed methods.

The basic idea of the post-processing formulations is to use the
optimal stability and convergence properties of the classical
Galerkin approximation for the potential field to derive  more
accurate velocity approximations than that obtained by direct
application of Darcy's law. In \cite{CORDES92} a post-processed
continuous flux (normal component of the velocity) distribution over
the porous media domain is derived using patches with
flux-conserving boundaries. Another conservative post-processing
based on a control volume finite element formulation is presented in
\cite{DURLOFSKY94}. Within the context of adaptive analysis,   a
local post-processing, is introduced in \cite{ZIENK92}, named
Superconvergent Patch Recovery, to recover higher order
approximations for the gradients of finite element solutions. %%

Mixed methods are based on the simultaneous approximation of
potential and velocity fields, and their main characteristic is the
use of different spaces for velocity and potential to satisfy a
compatibility  between the finite element spaces, which reduces the
flexibility in constructing finite approximations (LBB condition -
\cite{BREZZI74,FORTIN91}). A well known successful approach is the
dual mixed formulation developed by Raviart and Thomas
\cite{RAVIART77} using divergence based finite element spaces for
the velocity field combined with discontinuous Lagrangian spaces for
the potential. Stabilized mixed finite element methods have been
proposed to overcome the compatibilty conditions between the finite
element spaces, see \cite{BREZZI2005,LOULA88,MASUD2002} and
references therein.
In general, stabilized formulations use Lagrangian finite element
spaces and have been successfully employed in simulating Darcy's
flows in homogeneous porous media. The point is their applicability
to problems where the porous formation is composed by subdomains
with different conductivities. On the interface between these
subdomains, the normal component of Darcy velocity must be
continuous (mass conservation) but the tangential component is
discontinuous, and {any} formulation based on $C^0$ Lagrangian
interpolation for velocity fails  in representing the tangential
discontinuity, producing inaccurate approximations and spurious
oscillations, while completely discontinuous approximations for
velocity do not necessarily fit the continuity of the flux.
Interesting alternatives are the discontinuous Galerkin methods
\cite{BREZZI2005,HUGHES2005}, or coupling continuous and discontinuous Galerkin
methods as in \cite{DAWSON2005}.

In this work we will be concened with post-processing techniques derived form
variational formulations with good stability and convergence
properties on La\-gran\-gian finite element spaces.
In particular, we consider those based on stabilized finite element methods, such as the
global $C^0$ post-processing derived in \cite{ELSON90} by adding a least-square
residual of the mass conservation to a weak form of Darcy's law, and
the local post-processing (by element or macroelement) introduced in \cite{LOULA95}
which combines least squares residuals of both conservation of mass and
irrotationality condition with a discrete least squares residual of
Darcy's law taken at the superconvergence points of the gradient of the potential.
These formulations give velocity approximations
with higher convergence orders than  those obtained by the simple
direct use of Darcy's law. One important feature of the global
post-processing approach is the flexibility in the choice of
finite element spaces, allowing the use of equal-order Lagrangian
spaces, for example. In \cite{LOULA99} miscible displacement
analyses are performed using $C^0$ Lagrangian interpolations of the
same order for potential, velocity and concentration.

Here we study the approximation of flows in porous media
composed by layers with different conductivities, using
Lagrangian-based stabilized finite element methods.
We show that the Global Post-Processing presented in
\cite{ELSON90,LOULA95} can be derived from a Galerkin Least-Squares
stabilized mixed formulation and  popose generalized versions of the
global and local post-processing in which equal-order Lagrangian
finite element spaces are naturally adopted to porous media with
interfaces of discontinuity of the hydraulic conductivity.
For smooth interfaces of discontinuity, the proposed generalization
of the global $C^0$ post-processing  preserves the continuity of the
flux and exactly imposes the constraint between the tangent
components of Darcy velocity on the interface of the layers.
By construction, the local post-processing naturally
accommodates discontinuities of the tangent component but it does
not ensure  continuity  of the normal component of the velocity
field on the macroelement interfaces. When the interface of
discontinuity  is put in the interior of a macroelement, the
continuity or discontinuity  constraints on the normal and tangential
components of Darcy's velocity  are  exactly imposed.

After introducing some notations and definitions in Section 2, we
present in Section 3 our model problem for the flow in a porous
media with a smooth interface of discontinuity of the hydraulic
conductivity, consisting in a system of partial differential
equation with appropriate boundary and interface conditions. The
standard Galerkin  finite element method is reviewed in Section 4.
Stabilized mixed finite element methods are briefly discussed in
Section 5. In Section 6 we associate the origin of the considered
post-processing techniques with stabilized mixed formulations and
propose a generalization of  global and local
post-processings   to porous media with interface of discontinuity.
Numerical experiments illustrating the performance of the
global and local post-processing  are reported in Section 7 and in Section 8 we
draw some conclusions.

%
%-----------------------------------------------------------------------------
%

\section{Notations and Definitions}
\label{sec:notations}

Given $\Omega$ a bounded subdomain of $ \mathbb{R}^2$, with a sufficiently
smooth boundary $\partial \Omega$, let $L^2(\Omega)$ be the space of square
integrable functions in $\Omega$. For $m\geq 0$ integer, we define
$H^m(\Omega)$ the Hilbert space over $\Omega$ of order $m$
\[
H^{m}=\{v \in L^2(\Omega); D^{\alpha}v \in L^2(\Omega), |\alpha|\leq m\}
\mbox{, where}
\]
\[
D^\alpha(\cdot):=\frac{\partial^{|\alpha|}(\cdot)}
{\partial_{x_1}^{\alpha_1} \partial_{x_2}^{\alpha_2}},\quad
|\alpha| = \alpha_1+\alpha_2, \; \mbox{ and } \;
\alpha=(\alpha_1,\alpha_2) \mbox{ is an array}
\in \mathbb{N}^2.
\]
The inner product and the norm in $H^m(\Omega)$ are given by
\[
(v,w)_{m}:= \sum_{|\alpha| \leq m}
 \int_{\Omega} D^\alpha v  \, D^\alpha w \dx,
\quad \| v \|_{m}=(v,v)_{m}^{1/2}.
\]
We denote
\(
H^1_0(\Omega)=\{ v \in H^1(\Omega);\left.
v\right|_{\partial \Omega}=0 \}
\)
and define the space
\[
{H(\div)}= \left\{\bm{u} \in \left( L^2(\Omega) \right)^2;
\; \Div \bm{u} \in L^2(\Omega)\right\}, \mbox{ with norm }
\]
\[
\|\bm{u} \|_{H(\div)} := \left\{\| \bm{u} \|^2 + \| \Div \bm{u} \|^2\right \}^{1/2} .
\]
Finally, $H^0(\Omega)=L^2(\Omega)$ has
norm  $\|\cdot\|_0=\|\cdot\|$ and inner product given by
\[
(v,w) := \int_\Omega v\,w \dx.
\]

%
%-----------------------------------------------------------------------------
%

\section{Model Problem}
\label{sec:model}

\begin{figure}[htb]
\centering
\includegraphics[width=.5\linewidth]{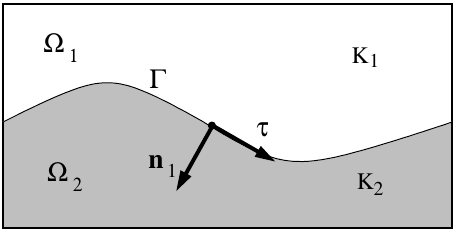}
\caption{Domain of the model problem.}
\label{fig:meios}%
\end{figure}

Let $\Omega = \Omega_1 \cup \Omega_2 \subset \mathbb{R}^2$, with
smooth boundary $\partial \Omega$ and outward unit normal vector
$\bm{n}$, be the domain of a rigid porous media saturated with an
incompressible homogeneous fluid. The domain $\Omega$ is composed by
two subdomains $\Omega_1$ and $\Omega_2$ that represent regions with
different conductivities, jointed by a smooth interface $\Gamma$, as
illustrated in Fig \ref{fig:meios}.

The mass balance asserts that
\begin{equation}
  \label{eq:balanco}
  \Div \bm{u} =f \quad \mbox{ in } \Omega
\end{equation}
where $f$ is a source of fluid and $\bm{u}$ is an average velocity of the fluid
in the porous media, given by Darcy's law
\begin{equation}
  \label{eq:darcy}
  \bm{u}=-K\nabla p \quad \mbox{ in } \Omega
\end{equation}
that relates the velocity field with the gradient of a hydraulic
potential $p$ through the hydraulic conductivity tensor $K$. For
simplicity, we consider a piecewise homogeneous porous media with
conductivities given by {$K=K_1$} in $\Omega_1$ and  {$K=K_2$} in
$\Omega_2$.
Restricting equations (\ref{eq:balanco}) and (\ref{eq:darcy}) to the
subdomains $\Omega_1$ and $\Omega_2$,  denoting
$\bm{u}_i=\bm{u}|_{\Omega_i}$, $p_i=p|_{\Omega_i}$, $f_i=f|_{\Omega_i}$,
$\Gamma_i:=\partial\Omega_i \backslash \Gamma$ ($i=1,2$) and
considering homogeneous Dirichlet boundary conditions for potential,
we have the following system of differential equations
\begin{eqnarray}
%\begin{array}{ll}
\bm{u}_1=-K_1 \nabla p_1  & \mbox{ in }  \quad \Omega_1
\label{eq:darcy1} \\
\bm{u}_2=-K_2 \nabla p_2  & \mbox{ in }  \quad \Omega_2
\label{eq:darcy2} \\
\Div  \bm{u}_1 =f_1           & \mbox{ in }  \quad \Omega_1
\label{eq:cm1}\\
\Div  \bm{u}_2 =f_2          & \mbox{ in }  \quad \Omega_2
\label{eq:cm2}
\end{eqnarray}
{with boundary conditions}
\begin{eqnarray}
p_1=0                    & \mbox{ on }  & \Gamma_1 \\
p_2=0                     & \mbox{ on } & \Gamma_2
\end{eqnarray}
{and interface conditions}
\begin{eqnarray}
K_1\nabla p_1\cdot \bm{n}_1 + K_2\nabla p_2 \cdot \bm{n}_2 = 0
&  \mbox{ on } & \Gamma
\label{eq:fluxnorm} \\
p_1=p_2                    & \mbox{ on } &  \Gamma.
\label{eq:cont}
\end{eqnarray}

Equation  (\ref{eq:fluxnorm}) establishes that the mass efflux in any
point on the interface vanishes, that is, the normal components of
the velocity field on $\Gamma$ must be continuous. Thus, we can
write
\begin{equation}
  \label{eq:fluxn}
  \bm{u}_1 \cdot \bm{n} = \bm{u}_2 \cdot \bm{n} \quad  \mbox{ on }  \quad
  \Gamma\ .
\end{equation}
Equation (\ref{eq:cont}) states the continuity of the hydraulic
potential. Multiplying (\ref{eq:darcy1}) and (\ref{eq:darcy2}) by
the unit vector $\tau$, tangent to $\Gamma$, we have:
\[
K_1^{-1}\bm{u}_1\cdot \tau =
-\nabla p_1 \cdot \tau = - \frac{\partial p_1}{\partial\tau},
\]
\[
K_2^{-1}\bm{u}_2\cdot \tau =
-\nabla p_2 \cdot \tau = - \frac{\partial p_2}{\partial\tau},
\]
and using equation  (\ref{eq:cont}) yields
\begin{equation}
\label{eq:fluxt}
K_1^{-1}\bm{u}_1\cdot \tau = K_2^{-1}\bm{u}_2\cdot \tau
\end{equation}
which shows that the tangent component of the velocity field should
be discontinuous on the interface $\Gamma$ for $K_1 \ne K_2$.

%
%--------------------------
%

\section{Galerkin Finite Element Approximations}
\label{sec:aprox}

As we will compute velocity approximations through post-processing
techniques that use {\it a priori} approximations of the potential
field, we first consider the single field problem posed by
substituting Darcy's law (\ref{eq:darcy}) in the mass conservation
(\ref{eq:balanco}):

\vspace{.2cm}
\paragraph{Problem P} {\it
Given the hydraulic conductivity  $K$ and a function $f$,
find the potential field $p$ such that
  \begin{equation}  \label{eq:p}
    -\Div \left(K \nabla p\right)   = f \quad  \mbox{ in }  \Omega
  \end{equation}
  with boundary condition
  \begin{equation}
    p  = 0 \quad \mbox{ on }  \partial \Omega.
  \end{equation}
}
Multiplying  (\ref{eq:p}) by a function $q \in H^1_0(\Omega)$ and
integrating it by parts over $\Omega$ we have the classical
variational formulation of  { Problem P}:

\paragraph{Problem PV} {\it
Find $p \in \mathcal{Q}=H^1_0(\Omega)$ such that
\begin{equation}
a(p,q)= f(q)
\quad \forall   q \in \mathcal{Q}
\end{equation}
with
\begin{equation}
a(p,q)=  \left(K \nabla p,\nabla q\right)
\, ; \quad
f(q)= (f,q).
\end{equation}  %
}
Letting $f\in H^{-1}(\Omega)$ we have the continuity of $f(\cdot)$.
As $a(\cdot,\cdot)$ is a continuous and $ \mathcal{ Q}-$elliptic
bilinear form, the problem has a unique solution by Lax-Milgram
Lemma.

\subsection{Galerkin Method for the Potential}
Let  $\{ \mathcal{T}_h \}$ be a family of partitions $\mathcal{
T}_h=\{\Omega^e\}$ of $\Omega$ indexed by the parameter $h$
representing the maximum diameter of the elements $\Omega^e \in \mathcal{T}_h$. Defining $ \mathcal{Q}_h^k \subset \mathcal{Q}$ as
\[
\mathcal{Q}_h^k=\{ q_h \in C^0(\Omega); \left. q_h\right|_{\Omega^e} \in
P_k(\Omega^e) \} \cap H_0^1(\Omega),
\]
the $C^0$ Lagrangian finite element space of degree  $k\geq 1$ in
each element $\Omega^e$, where $P_k(\Omega^e)$ is the set of the
polynomials of degree $\leq k$ posed on $\Omega^e$, we present the
Galerkin finite element approximation  to the {Problem PV}
\paragraph{Problem Ph} 
{\it
Find $p_h \in \mathcal{Q}_h^k$ such that
\begin{equation}
\label{eq:pvh}
a(p_h,q_h)= f(q_h)
\quad \forall   q_h \in \mathcal{Q}_h^k.
\end{equation}
}

As $ \mathcal{Q}_h^k \subset \mathcal{Q}$,  {Problem Ph} has
existence and uniqueness of  solution guaranteed and the following
estimate holds (see, for example, in Ciarlet \cite{CIARLET78})
\begin{equation}
  \label{eq:errop}
  \| p - p_h \| + h \| \nabla p - \nabla p_h \| \leq C h^{k+1}
  \left| p \right|_{k+1} \, .
\end{equation}

%
%----------------------------------
%
\subsection{Galerkin Method for the Velocity}

After solving { Problem Ph} we can evaluate the velocity field
directly through Darcy's law (\ref{eq:darcy})
\begin{equation}
  \label{eq:vg}
  \bm{u}_G := -K \nabla p_h.
\end{equation}

From (\ref{eq:errop}) we have the following error bound to this approximation
\begin{equation}
  \label{eq:errovg}
  \| \bm{u} - \bm{u}_G \| \leq  Ch^k | p|_{k+1}.
\end{equation}

This direct approximation for the velocity field is, in principle,
completely discontinuous on the interface of the elements and
presents very poor mass conservation and a too weak approximation of
flux (Neumann) boundary conditions. In \cite{CORDES92,DURLOFSKY94},
locally conservative post-processings are presented that use
(\ref{eq:vg}) to construct continuous flux fields over the domain.
Totally continuous velocity approximations can be obtained using the
classical post-precessing approach based on a continuous $L^2$
smoothing of $K\nabla p_h$. For continuous velocity fields this
simple post-processing is very accurate with linear or bilinear
elements. For higher order element it may be less accurate than the
direct approximation.

%
%--------------------------
%
\section{Stabilized Mixed Formulations}
\label{sec:stbf}

Before focusing on the post-processing techniques, we present the
basic concepts of consistenly stabilized mixed methods.
In its velocity and potential formulation Darcy problem, given by equations
(\ref{eq:balanco}) and (\ref{eq:darcy}), fits in the abstract mixed
formulation studied by Brezzi \cite{BREZZI74}
\[
A \bm{u} + B^*p = g \ \ \ {\rm in } \ \ \Omega
\]
\[
B \bm{u} = f \ \ \ {\rm in } \ \ \Omega
\]
whose weak form consists in:
Find $\{\bm{u},p \} \in W\times V$ such that
\begin{equation}
a(\bm{u},\bm{v}) + b(\bm{v},p ) - g(\bm{v}) + b(\bm{u},q)-f(q) =0 \ \ \forall \;
\{\bm{v},q\} \in W\times V
\end{equation}
with
\[
a(\bm{u},\bm{v}) = (A\bm{u},\bm{v})  \ \ \forall \  \bm{v} \in W \, ,
\]
\[
b(\bm{v},q ) = (B\bm{v},q) \ \ \forall \ \{\bm{v},q\} \in W\times V \, ,
\]
\[
g(\bm{v}) =(g,\bm{v}) \ \ \forall \  \bm{v}\in W \, ,  \quad
f(q) =(f,q) \ \ \forall \  q\in V
\]
where $W$ and $V$ are appropriate Hilbert spaces. For continuous
linear and bilinear forms, existence and uniqueness of solutions
of this abstract mixed formulation are assured by the following
well known Babu\v ska - Brezzi or  inf-sup conditions:
\begin{equation}
\sup_{\bm{v}\in W} \frac{a(\bm{u},\bm{v}) }{ \| \bm{v} \|_W } \ge \alpha \| \bm{u}
\|_W \ \ \forall \ \bm{v} \in K_0 \, ,
\end{equation}
\begin{equation}
\sup_{\bm{v}\in W} \frac{b(\bm{v},q) }{ \| \bm{v}\|_W } \ge \beta \| q\|_V \, ,
\ \ \forall \; q \in V \, ,
\end{equation}
with
\begin{equation}
K_0= \{ \bm{v}\in W \, , \ b(\bm{v},q ) = 0 \ \ \forall\; q \in V\} \, .
\end{equation}
These kind of compatibility conditions usually impose severe
limitations in the choice of stable finite element approximations.
To overcome these limitations stabilized finite element
formulations have been proposed. Here we will consider
the Galerkin Least-Squares (GLS) or Mixed Petrov-Galerkin method
introduced in  \cite{LOULA87,FRANCA88} consisting in:
{\it
Find $\{\bm{u}_h,p_h \} \in W_h\times V_h \subset W\times V$ such that
$\forall\;\{\bm{v}_h,q_h\} \in W_h\times V_h$
\[
a(\bm{u}_h,\bm{v}_h) + b(\bm{v}_h,p_h ) - g(\bm{v}_h) + b(\bm{u}_h,q_h)-f(q_h) +
\]
\[
{\delta_1(A \bm{u}_h + B^*p_h - g, A\bm{v}_h + B^* q_h)_h + \delta_2
(Bu_h-f,Bv_h)_h} = 0
\]
}
\noindent where $\delta_1$ and $\delta_2$ are real parameters activating the
least squares residuals of the governing equations in the interior
of the elements.

\subsection{ Turning Saddle Point into a Minimization Problem }

Let
\begin{equation}
\mathcal{U}={H(\div)} \  \hbox{ and} \
 \mathcal{Q}=H^1_0(\Omega) \, ,
\end{equation}
and $ \mathcal{U}_h^l \subset \mathcal{U}$ and $\mathcal{Q}_h^k \subset \mathcal{Q}$ typical $C^0$
Lagrangian finite element spaces of degree $l$ and $k$,
respectively. Within these spaces the following stabilization is
proposed in \cite{LOULA88} for a heat transfer problem with the same
mathematical structure of Darcy's problem
\paragraph{Problem GLS}
{\it
Find $\{\bm{u}_h,p_h\} \in \mathcal{U}_h^l\times\mathcal{Q}_h^k$ such that
\[
B_\mathrm{GLS}(\{\bm{u}_h,p_h\},\{\bm{v}_h,q_h\}) =
     \delta_2(\lambda f,\Div \bm{v}_h)- (f,q_h)
       \quad  \forall\;\left\{\bm{v}_h,q_h \right\} \in \mathcal{U}_h^l\times \mathcal{Q}_h^k
\]
}
with
\begin{eqnarray*}
  B_\mathrm{GLS}(\{\bm{u},p\},\{\bm{v},q\}) & = & (\lambda \bm{u},\bm{v}) - (\Div \bm{v}, p)
  - (\Div \bm{u}, q) \\
  & & +\delta_1 \left( K\left( \lambda \bm{u} + \nabla p\right),
    \lambda \bm{v} + \nabla q\right)
   +\delta_2 (\lambda \Div \bm{u}, \Div \bm{v})  \, .
\end{eqnarray*}

For appropriate choices of $\delta_1>0$ and  $\delta_2>0$, {
Problem GLS} is equivalent to the minimization problem:
Find $\{\bm{u}_h,p_h\} \in \mathcal{U}_h^l\times\mathcal{Q}_h^k$ such that
\[
 J(\{\bm{u}_h,p_h\}) \le  J(\{\bm{v}_h,q_h\})\ \ \forall\  \{ \bm{v}_h,q_h\} \ \in \mathcal{U}_h^l\times \mathcal{Q}_h^k
\]
with
\[
 J(\{\bm{v}_h,q_h\})= \frac{1}{2}B_\mathrm{GLS}(\{\bm{v}_h,q_h\},\{\bm{v}_h,q_h\}) +
 (f,q_h) - \delta_2( \lambda f,\Div \bm{v}_h)  \, ,
\]
which guarantees  existence, uniqueness and the best approximation
property for this stabilized formulation in the energy norm. The
error bound for this method is given by
\begin{eqnarray}
\nonumber
\| \bm{u} -\bm{u}_h\|_{H(\div)} + \|p-p_h\|_1 \leq C \left(
 \| \bm{u} -\bm{v}_h\|_{H(\div)} + \|p-q_h\|_1 \right) \\
 \forall \; \{\bm{v}_h,q_h\} \in
 \mathcal{U}_h\times\mathcal{Q}_h.
\label{eq:boundgls}
\end{eqnarray}
For sufficiently regular exact solutions, same order $C^0$
Lagrangian spaces ($l=k$) lead to the error estimates
\begin{equation}
    \label{eq:mpe-p}
     \|p - p_h\| +h\| \nabla p - \nabla p_h  \|\leq  C
     h^{k+1}\left( \left|\bm{u} \right|_{k+1}
     + \left|p \right|_{k+1} \right),
\end{equation}
\begin{equation}
    \label{eq:mpe-div}
     \|\bm{u} - \bm{u}_h \|_{{H(\div)}} \leq  C  h^k\left( \left|\bm{u} \right|_{k+1}
     + \left|p \right|_{k+1} \right),
\end{equation}
\begin{equation}
    \label{eq:mpe-u}
     \|\bm{u} - \bm{u}_h \| \leq  C  h^k\left( \left|\bm{u} \right|_{k+1}
     + \left|p \right|_{k+1} \right),
\end{equation}
whose orders of convergence are optimal for potential in
$L^2(\Omega)$ norm and $H^1(\Omega)$ seminorms (equivalent to
$H^1(\Omega)$ norm), and for velocity in $ H(\hbox{div}) $ norm, but
suboptimal for the velocity in $[L^2(\Omega)]^2$ norm.

\subsection{Adjoint Formulation}

More recently Masud and Hughes introduced in \cite{MASUD2002} the
following non symmetric stabilization  for Darcy flow
\vspace{.2cm}

\paragraph{Problem HVM} 
{\it
       Find $\{\bm{u}_h,p_h\} \in \mathcal{U}_h^l\times\mathcal{Q}_h^k$ such that
    \[
      B_\mathrm{HVM}(\{\bm{u}_h,p_h\},\{\bm{v}_h,q_h\}) =  (f,q_h)
       \quad  \forall\;\left\{\bm{v}_h,q_h \right\} \in \mathcal{U}_h^l\times \mathcal{Q}_h^k
    \]
      }
with
\begin{eqnarray*}
   B_\mathrm{HVM}(\{\bm{u},p\},\{\bm{v},q\}) &=& (\lambda \bm{u},\bm{v}) - (\Div \bm{v}, p) + (\Div \bm{u}, q)
\\
   &&+\frac{1}{2} \left( K\left( \lambda \bm{u} + \nabla p\right),
  -\lambda \bm{v} + \nabla q\right) \, .
\end{eqnarray*}

According to the authors of  \cite{MASUD2002}  stabilizations such
as Galerkin/Least Squares (GLS) are not as effective for the current
problem as this adjoint formulation. The only essential difference
is the sign on the $\bm{v}$ term which is considered crucial. In fact
due to the non symmetry one can easily prove stability of this
formulation in the sense of Lax, since
\[
   B_\mathrm{HVM}(\{\bm{v},q\},\{\bm{v},q\}) \ge \frac{\lambda}{2} \left\|\bm{v}\right\|^2
   + \frac{K}{2} \left\| \nabla q \right\|^2
 \ge \alpha (\| \bm{v} \|^2 + \| \nabla q \|^2) \ \ \forall \left\{\bm{v},q \right\} \in
\mathcal{U} \times \mathcal{Q}
\]
with $\alpha = \min\{ \lambda , K \}/2$,  for  $K$ and $\lambda$ positive constants.
Thus, the HVM formulation is unconditionally stable in $[L^2(\Omega)]^2 $
norm for
velocity and $H^1(\Omega)$ seminorm for potential, while the
GLS formulation introduced before is stable in ${H(\div)}$ norm and
$H^1(\Omega)$ seminorm for velocity and potential, respectively.
Choosing $\delta_1 =-1/2$, $ \delta_2=0$ in the GLS formulation,
or taking $q_h= -q_h$ in HVM formulation,
restricted to the product
space $ \mathcal{U}_h^l\times\mathcal{Q}_h^k \subset \bm{u} \times \mathcal{Q}$, with
\[
    \bm{u}=[L^2(\Omega)]^2 \ , \  \mathcal{Q}=H^1_0(\Omega)  \, , \quad
\| \left\{\bm{v},q \right\} {\|^2_{_{\bar{\mathcal{U}} \times {\mathcal{Q}}}}} = \| \bm{v} \|^2 + \| \nabla q \|^2 \, ,
\]
and  using integration by parts, both GLS and HVM methods reduce
to the following Galerkin Least Squares formulation
\vspace{.2cm}
\paragraph{Problem GLS1} 
{\it
       Find $\{\bm{u}_h,p_h\} \in \mathcal{U}_h^l\times\mathcal{Q}_h^k$ such that
    \[
      B_\mathrm{GLS1}(\{\bm{u}_h,p_h\},\{\bm{v}_h,q_h\}) = -(f,q_h)
     {  \quad  \forall\;\left\{\bm{v}_h,q_h \right\} \in \mathcal{U}_h^l\times \mathcal{Q}_h^k}
    \]
      }
where
\begin{eqnarray*}
B_\mathrm{GLS1}(\{\bm{u},p\},\{\bm{v},q\}) &=& (\lambda \bm{u},\bm{v}) + (\bm{v},\nabla p) + (\bm{u},\nabla q) 
\\
&&-\frac{1}{2} \left( K\left( \lambda \bm{u} + \nabla p\right),
           \lambda \bm{v} + \nabla q\right)
\end{eqnarray*}
is symmetric but not stable in the sense of Lax.
However, this particular GLS method preserves stability in the
sense of Babu\v ska, since we can always choose
\[
\left\{ \bar{\bm{v}}, \bar{q} \right\} =
    \left\{ \bm{u}, -p \right\}
\]
such that  $\forall \; \left\{\bm{u},p \right\} \in \bm{u}\times \mathcal{Q}$
   \[
    \sup_{\left\{\bm{v},q \right\}\in \bm{u}Q}\frac{|B_\mathrm{GLS1}(\left\{\bm{u},p \right\},\left\{\bm{v},q \right\})|}{\|\left\{\bm{v},q \right\}{\|_{_{\bar{\mathcal{U}} \times {\mathcal{Q}}}}}}\geq
    \frac{|B_\mathrm{GLS1}(\left\{\bm{u},p \right\},\left\{ \bar{\bm{v}}, \bar{q} \right\}|}
    {\|\left\{ \bar{\bm{v}}, \bar{q} \right\}{\|_{_{\bar{\mathcal{U}} \times {\mathcal{Q}}}}}} =
     \frac{ \lambda \|\bm{u}\|^2 + K \|\nabla p \|^2} {2 \|
     \{\bm{u},p\}{\|_{_{\bar{\mathcal{U}} \times {\mathcal{Q}}}}} } \, ,
    \]
yielding
   \[
    \sup_{\left\{\bm{v},q \right\}\in \bm{u}Q}\frac{|B_\mathrm{GLS1}(\left\{\bm{u},p \right\},\left\{\bm{v},q \right\})|}{\|\left\{\bm{v},q \right\}{\|_{_{\bar{\mathcal{U}} \times {\mathcal{Q}}}}}}\geq
     \alpha \| \left\{\bm{u},p \right\}   {\|_{_{\bar{\mathcal{U}} \times {\mathcal{Q}}}}} \, ,
    \]
with  again $\alpha = \min\{ \lambda , K \}/2$, as in the HVM formulation.
 We should observe that the above stability result is
proved independently of any compatibility condition between the
spaces $\bm{u}$ and $\mathcal{Q}$. Therefore, GLS1 is consistent and
unconditionally stable in the sense of Babu\v ska which means that
any conforming GLS1 finite element approximation is stable.  Now,
the error bound is given by
\begin{equation}
\label{eq:boundhvm} \| \bm{u} -\bm{u}_h\| + \|p-p_h\|_1 \leq C \left(
 \| \bm{u} -\bm{v}_h\| + \|p-q_h\|_1 \right)\quad \forall \; \{\bm{v}_h,q_h\} \in
 \mathcal{U}_h \times\mathcal{Q}_h.
\end{equation}
For
sufficiently regular exact solutions, same order $C^0$ Lagrangian
spaces ($l=k$) lead to the error estimates
\begin{equation}
\label{eq:hvm-p}
     \|p - p_h\| +h\| \nabla p - \nabla p_h  \|\leq  C
     h^{k+1}\left( \left|\bm{u} \right|_{k+1}
     + \left|p \right|_{k+1} \right),
\end{equation}
\begin{equation}
\label{eq:hvm-div}
     \|\bm{u} - \bm{u}_h \|_{{H(\div)}} \leq  C  h^{k-1}\left( \left|\bm{u} \right|_{k+1}
     + \left|p \right|_{k+1} \right),
\end{equation}
\begin{equation}
\label{eq:hvm-u}
     \|\bm{u} - \bm{u}_h \| \leq  C  h^k\left( \left|\bm{u} \right|_{k+1}
     + \left|p \right|_{k+1} \right),
\end{equation}
with  optimal orders for potential in $L^2(\Omega)$ and
$H^1(\Omega)$ norms but suboptimal for  velocity  in
$[L^2(\Omega)]^2$ and $H(\hbox{div})$ norms. The same orders of convergence are obtained
for the non symmetric HVM formulation as a consequence of Lax and
C\'ea's lemmas.

%
%-------------------------------------------------------------
%

\section{Continuous/Discontinuous Post-Processings}
\label{sec:cdpp}

We start the presentation of the post-processing formulations, by identifying
their origin with stabilized mixed finite element methods. Observing
that the least-squares residual of Darcy's law in the GLS
formulation is responsible for the stabilization of the potential in
$H^1$ seminorm, we consider the following GLS stabilization: {\it
Find $\{\bm{u}_h,p_h\} \in \mathcal{U}_h^l\times\mathcal{Q}_h^k$ such that
\begin{eqnarray*}
(\lambda \bm{u}_h,\bm{v}_h) &-& (\Div \bm{v}_h, p_h) - (\Div \bm{u}_h, q_h)
+(f,q_h)
\\
&+&\delta_1\left[\left( K\nabla p_h , \nabla q_h\right) -
(f,q_h)\right] +\delta_2 (\lambda(\Div \bm{u}_h -f), \Div \bm{v}_h) 
=0
\end{eqnarray*}
}
in which the least squares residual of Darcy's law is replaced
by the Galerkin residual of
the potential equation. Again we have derived a consistent
and unconditionally stable  formulation
with velocity in $H(\hbox{div})$ and potential in $H^1(\Omega)$.
Alternatively, we can solve first
the potential equation independently of the velocity,
using the classical Galerkin method, and then
compute the velocity using the GLS method with the given
approximation for the potential.
This is exactly the strategy adopted in the post-processing
techniques introduced in
\cite{ELSON90} and analyzed in \cite{LOULA95}, as  shown next.
\subsection{Global $C^0$ Post-Processings}
\label{sec:gpp}
Let  $\mathbb{W}_h$ be a finite element subspace of ${H(\div)}$ and
$p_h$ be the solution of {Problem Ph}.
The global post-processing is based on the following residual form
\[
\left(K^{-1}\bm{u}_h + \nabla p_h, \bm{w}_h \right) +
(\delta h)^\alpha\left( \Div \bm{u}_h - f, \Div \bm{w}_h \right) = 0 \quad
\forall \; \bm{w}_h \in\mathbb{W}_h,
\]
with  $\delta$ and $\alpha$ real positive parameters. For fixed
$\delta$ we introduce the bilinear form
\[
b_\alpha(\bm{u}_h,\bm{w}_h)=
\left(K^{-1}\bm{u}_h ,\bm{w}_h \right) +
(\delta h)^\alpha\left( \Div \bm{u}_h, \Div \bm{w}_h \right)
\]
and define a family of Global Post-Processings (GPP), depending on
$\alpha$, as follows:

\vspace{.2cm}
\paragraph{Problem GPP}{\it
  Given $p_h$ an approximation to the potential field,
  find the post-processed velocity field
  $\bm{u}_h\in\mathbb{W}_h $ such that
  \begin{equation}
    \label{eq:pp}
    b_\alpha(\bm{u}_h,\bm{w}_h)= (\delta h)^\alpha\left( f, \Div \bm{w}_h \right) -
    \left(\nabla p_h, \bm{w}_h \right)
    \quad \forall  \; \bm{w}_h   \in\mathbb{W}_h.
  \end{equation}
}
The bilinear form $b_\alpha(\cdot,\cdot)$ is continuous, symmetric
and positive definite. This fact, associated with conforming
approximations, guarantees existence and uniqueness of the solution
for {Problem GPP}. A complete analysis of this method can be
found in \cite{LOULA95} where the following estimate is proved
\begin{eqnarray}
  \| \bm{u}-\bm{u}_h\| + (\delta h)^{\alpha/2}\| \Div \bm{u} - \Div \bm{u}_h \| &\leq &
  C \left(\|\bm{u}-\bm{w}_h\| +(\delta h)^{\alpha/2}\| \Div \bm{u} - \Div \bm{w}_h \|
  \right. \nonumber \\
  & &
  \left. +(\delta h)^{-\alpha/2}\| p - p_h \| \right) \quad \forall  \; \bm{w}_h
  \in\mathbb{W}_h,
  \label{eq:est1}
\end{eqnarray}
with $C$ independent of $h$.
The application of this formulation associated with $C^0$ Lagrangian
subspaces of $H^1$:
\[
\mathcal{W}_h^l=
\{\varphi_h \in C^0(\Omega); \; \left. \varphi_h\right|_{\Omega^e} \in
P_l(\Omega^e)
\; \forall \;\Omega^e \in \mathcal{T}_h \} \, ,
\]
\[
\mathbb{W}_h^l= \{ \bm{w}_h \in \mathcal{W}_h^l \times \mathcal{W}_h^l \} \, ,
\]
leads to the error estimate
\begin{equation}
  \| \bm{u}-\bm{u}_h\| + (\delta h)^{\alpha/2}\| \Div \bm{u} - \Div \bm{u}_h \|\leq
  C \left(h^{l+\alpha/2}|\bm{u}|_{l+1} + h^{k+1-\alpha/2}|p|_{k+1} \right)
  \label{eq:est2}
\end{equation}
for sufficiently regular exact solution. Of course, this choice
gives completely continuous velocity approximations incompatible
with heterogeneous media with material discontinuities.  In next
section we present  global and local post-processing formulations in
which the velocity field  can be continuous or discontinuous on the
interfaces of material discontinuities.

\subsection{ Global Post-Processings with Interface of Discontinuity}
\label{sec:ppge}

Using equations (\ref{eq:darcy1}-\ref{eq:cm2}), (\ref{eq:fluxn}) and
(\ref{eq:fluxt}) we can proceed as in the derivation of (\ref{eq:pp}) and
write the problem: {\it find $\bm{u}_h^1 \in \mathbb{W}_{h1}^l$ and
$\bm{u}_h^2 \in \mathbb{W}_{h2}^l$ such that
\[
b^1_\alpha(\bm{u}^1_h,\bm{w}_h)= (\delta h)^\alpha\left( f, \Div \bm{w}_h \right) -
\left(\nabla p^1_h, \bm{w}_h \right)
\quad \forall  \; \bm{w}_h   \in\mathbb{W}^l_h
\]
\[
b^2_\alpha(\bm{u}^2_h,\bm{w}_h)= (\delta h)^\alpha\left( f, \Div \bm{w}_h \right) -
\left(\nabla p^2_h, \bm{w}_h \right)
\quad \forall  \; \bm{w}_h   \in\mathbb{W}_h^l
\]
plus the interface conditions
\begin{equation}
  \label{eq:restr}
  \bm{u}^1_h\cdot \bm{n} = \bm{u}^2_h\cdot \bm{n} \quad \mbox{ and } \quad
  K_1^{-1} \bm{u}^1_h\cdot \tau= K_2^{-1} \bm{u}^2_h\cdot \tau
  \quad \mbox{ on } \; \Gamma
\end{equation}
}
where
\[
b^j_\alpha(\bm{v}_h,\bm{w}_h)=(K^{-1}_j \bm{v}_h, \bm{w}_h) + (\delta h)^\alpha
(\Div \bm{v}_h , \Div \bm{w}_h),
\]
\[
\mathcal{W}_{hj}^l= \{\varphi_h \in C^0(\Omega_j); \; \left.
\varphi_h\right|_{\Omega^e} \in P_l(\Omega^e) \; \forall \;\Omega^e
\in \mathcal{T}_h \},
\]
\[
\mathbb{W}_{hj}^l=  \{\bm{w}_h \in \mathcal{W}_{hj}^l \times \mathcal{W}_{hj}^l\}.
\]

One possible form to impose the relations on the velocity field
would be write  $\bm{u}_h^2=F(\bm{u}_h^1)$,  $\bm{u}_h^1=F^{-1}(\bm{u}_h^2)$ and
add residuals of (\ref{eq:restr}) on the interface $\Gamma$, to the
variational formulation using Lagrange multiplier, for example.
Alternatively, we propose a simpler finite element method which
incorporates the discontinuities on the interface of the elements
leading to a discrete problem similar to that corresponding to the
global $C^0$ post-processing.

\subsubsection*{ Incorporating the Discontinuity on the Interface}\
\label{sec:idi}

\begin{figure}[htb]
  \centering
  \includegraphics[width=.85\linewidth]{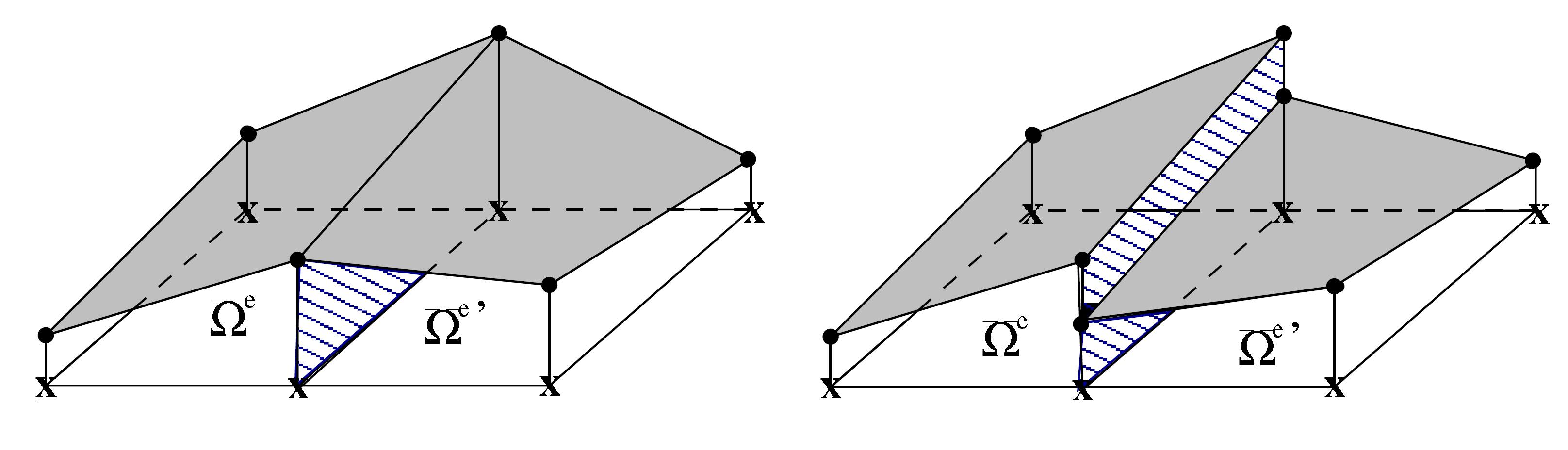}
  \caption{$C^0$ class elements and discontinuous ($C^{-1}$) elements.}
  \label{fig:elc}%
\end{figure}

Let $\bm{n}_1=\bm{n}$ be the outward normal to the boundary of $\Omega_1$
at a generic point  $P(\bm{x})$ on the interface $\Gamma$. Thus, for
$\epsilon>0$
\[
\bm{u}^1 = \lim_{\epsilon \rightarrow 0}\bm{u}(\bm{x}-\epsilon \bm{n})
\]
\[
\bm{u}^2 = \lim_{\epsilon \rightarrow 0}\bm{u}(\bm{x}+\epsilon \bm{n}) \, .
\]
According to Eq. (\ref{eq:fluxt})
\[
 \lim_{\epsilon \rightarrow 0} K^{-1}(\bm{x}-\epsilon \bm{n})\bm{u}(\bm{x}-\epsilon \bm{n})\cdot\tau=
K_1^{-1}\bm{u}^1 \cdot \tau = K_2^{-1}\bm{u}^2 \cdot \tau =
 \lim_{\epsilon \rightarrow 0} K^{-1}(\bm{x}+\epsilon \bm{n})\bm{u}(\bm{x}+\epsilon \bm{n})\cdot \tau
\]
which implies \( \displaystyle \bm{u}^1 \cdot \tau \ne \bm{u}^2 \cdot \tau \) for
$K_1 \ne K_2 $.  The $C^0$ Lagrangian interpolation for the velocity
field is obviously incompatible with this last condition since it
imposes exactly
\begin{equation}
  \label{eq:resc}
  \bm{u}^1_h\cdot \bm{n} = \bm{u}^2_h\cdot \bm{n} \quad \mbox{ and } \quad
  \bm{u}^1_h\cdot \tau= \bm{u}^2_h\cdot \tau
\end{equation}
on the interfaces of the elements, as illustrated in Fig
\ref{fig:elc}.
In this case, the approximate solution ${\bm{u}}_h$ at a global node
on $\Gamma$ will represent an intermediate value of the discontinuous
solution. Our idea is to impose exactly the constraint
\begin{equation}
  \label{eq:resh}
 K_1^{-1} \bm{u}^1_h\cdot \tau= K_2^{-1}\bm{u}^2_h\cdot \tau
\end{equation}
on the interface of discontinuity and assemble the global finite
element approximation in terms of a reference solution
$\bar{\bm{u}}_h$, which is uniquely defined at the nodes. This
procedure leads a finite element formulation with the same
connectivity of the global $C^0$ post-processing. We start by
choosing as reference solution ${\bm{u}}^2_h$, such that
\begin{equation}
\label{eq:uref1}
\bm{u}_h^2=\mathrm{I}\,\bar{\bm{u}}_h \quad \mbox{ on } \quad \Gamma\, ,
\end{equation}
where $\mathrm{I}$ is the identity.
We now relate $\bm{u}_h^1$ to the reference solution. From
Eqs (\ref{eq:restr}) and (\ref{eq:uref1}) we can write
\[
\mathrm{T}_1 \bm{u}_h^1= \mathrm{T}_2  \,\bm{u}_h^2 =  \mathrm{T}_2  \,\bar{\bm{u}}_h
\quad \mbox{ on } \quad \Gamma,
\]
with the Cartesian matrix of $ \mathrm{T}_i $ given by
\begin{equation}
\label{eq:tmatrix}
 \left[ \mathrm{T}_i \right] = \left[
\begin{array}{cc}
  \lambda_{11}^i \tau_1 + \lambda_{21}^i\tau_2  &
  \lambda_{12}^i \tau_1 + \lambda_{22}^i\tau_2  \\
   \mathrm{n}_1 & \mathrm{n}_2
\end{array}
\right]
\end{equation}
where $(\tau_1,\tau_2)$ and $(\mathrm{n}_1,\mathrm{n}_2)$ are the Cartesian components of
the unit vectors $\tau$ and $\bm{n}$, respectively, and $\lambda_{lm}^i$ are the
components of the hydraulic resistivity tensor $\lambda_i=K_i^{-1}$.
Considering that $\mathrm{T}_1$ is invertible, we get
\begin{equation}
\label{eq:uref2}
\bm{u}_h^1 =
\mathrm{T}_1^{-1}
\mathrm{T}_2\; \bar{\bm{u}}_h \quad \mbox{ on } \quad \Gamma.
\end{equation}
Using Eqs  (\ref{eq:uref1}) and  (\ref{eq:uref2}) we are able to
pose the interface problem in terms of a continuous global
vector of unknowns $\{\bar{u}\}$ only.  This can be done at element
level as follows: let $\{u^e\}$ be the array of unknowns of a
generic element $\Omega^e$ with just one edge belonging to
the interface of discontinuity, and $
\left\{ u^e_{\Gamma} \right\}$ the unknowns associated with the nodes
on that edge belonging to the interface.
By using (\ref{eq:uref2}) to each node on $\Gamma$, we
can construct a matrix $\left[ \displaystyle T_\Gamma \right]$ and write
\begin{equation}
  \label{eq:rel1}
  \displaystyle \left\{ \displaystyle u^{e^{\,}}_{\Gamma} \right\} = \left[ \displaystyle T_\Gamma \right]
  \{ \displaystyle u^{e'}_{\Gamma} \}
\end{equation}
where $\Omega^e$ and $\Omega^{e'}$ are elements that share the
interface of material discontinuity as illustrated in Fig
\ref{fig:elc}. Choosing $\{u^{e'}_{\Gamma}\}$ as reference we have
\[
%  \label{eq:uinc2}
  \left\{ u^e \right\} = \left\{
  \begin{array}{c}
    \left\{ u^e_{\bar{\Omega} \backslash \Gamma} \right\} \\
    --- \\
    \left[ \displaystyle T_\Gamma \right] \{ \displaystyle u^{e'}_{\Gamma} \}
  \end{array}
  \right\}
     = \left[ T \right]
  \left\{
  \begin{array}{c}
    \left\{ u^e_{\bar{\Omega} \backslash \Gamma} \right\} \\
    --- \\
    \{ \displaystyle u^{e'}_{\Gamma} \}
  \end{array}
  \right\}
\]
or
\begin{equation}
  \label{eq:trans1}
  \left\{ u^e \right\} = \left[ T \right] \{ \bar{u}^{e} \}
\end{equation}
where  $[T]$ is the matrix of the linear transformation from  $ \{
u^e \}$ to $\{ \bar{u}^{e} \}$. The contribution of the element
$\Omega^e$ to the global stiffness $[K]$ is given by:
\[
[K^e]\{u^e\}=[ K^e       ] [ T ] \{\bar{u}^e\}
            =[ \bar{K}^e ] \{\bar{u}^e\},  %=  \{f^e\}.
\]
with the new matrix  $[\bar{K}^e]$ non-symmetric.
The continuous/discontinuous approximation $\{ u^e \}$ on the interface
is obtained at element-level  using (\ref{eq:trans1}) after solving the global
system in  $\{\bar{u}\}$.

 \paragraph{\it Remarks.}
\begin{enumerate}

\item[(i)] In (\ref{eq:tmatrix}) we used unit vectors $\tau$ and $\bm{n}$
calculated at a nodal point on the interface. Using isoparametric
elements, the interface $\Gamma$  is approximated by a interpolant
$\Gamma_h \in C^0(\Gamma)$. In general, for curved interfaces these
vectors are not uniquely defined at the vertices of the elements. In
this case averages of these unit vectors can be adopted.

\item[(ii)] Considering that the solution of our model problem is piecewise
  smooth, we should expect that the error estimate (\ref{eq:est2})
  holds for each subdomain $\Omega_i$ when the  continuity/discontinuity constraints on
  the normal/tangential components of Darcy's velocity are exactly imposed on the interface
  of the two media with different hydraulic conductivity.
  A convergence study is presented in Sec. \ref{sec:ex3} confirming the orders of convergence
  predicted  in (\ref{eq:est2}).

\item[(iii)]  As the weak forms of the  global post-processings are posed
in finite element subspaces of ${H(\div)}$, they allow the use of
divergence based finite element spaces, such as Raviart-Thomas spaces,
naturally accommodating discontinuities of the tangent components of the
velocity field on the interface of the elements.

\item[(iv)] Although introduced for the Global Post-Processing formulation
    this idea can, of course, be applied to more general situations,
    as for example,  to the stabilized mixed finite element
    methods introduced in \cite{LOULA88,MASUD2002} or to the local post-processing introduced
    in \cite{LOULA95}.
\end{enumerate}

\subsection{ Local Post-Processings with Interface of Discontinuity}
A stabilized mixed formulation using totally discontinuous
Lagrangian interpolations for both velocity and potential fields is
presented and analyzed  in \cite{HUGHES2005, BREZZI2005}.
Stability is provided, for $k\ge 1$,
by consistently adding to an originally unstable DG formulation
least squares or adjoint residuals of Darcy's law. For using
discontinuos interpolation, the velocity degrees of freedom are
eliminated at element level in favor  of the potential degrees of
freedom. After solving  the global system corresponding to  the
potential,  the velocity approximation is computed by a local
(element by element) post-processing. As proved in
\cite{BREZZI2005}, choosing equal order interpolations this
stabilized mixed method leads to optimal orders of convergence for
the pressure and to the suboptimal estimate for the velocity
approximation
$$
 \| \bm{u} - \bm{u}_h \| \leq  C  h^{k} \| p \|_{k+1} \quad \mbox{ for} \quad k\ge 1 \, .
$$

Taking advantage of the superconvergence of the gradient of the
Galerkin finite element solution at special points in the interior
of the elements, a local
post-processing at element or macroelement level is presented in
\cite{LOULA95}, using least squares residuals of both  equilibrium
equation and the irrotationality condition. Let $ \mathcal{S}_h^k $ be
triangular or quadrilateral Lagrangian finite element subspaces of $
L^2(\Omega)$ consisting of piecewise polynomials of degree $k$ on
each element.
By construction,  we take $ \mathcal{S}_h^k $ of class $C^0$ in each macroelement but discontinuous
on the macroelement boundaries.  In the product space
$\mathbb{M}_{h}=\mathcal{ S}_h^k\times \mathcal{S}_h^k$  we consider the
mesh-dependent bilinear form
$$
 b_{h} (\bm{u}_h, \bm{v}_h) = ( \bm{u}_h, \bm{v}_h)_{h} +
 h^2 (\Div \bm{u}_h, \Div \bm{v}_h) +  h^2( \curl \bm{u}_h ,
\curl  \bm{v}_h) \quad \forall \bm{v}_h, \bm{u}_h \in \mathbb{M}_{h} \, ,
$$
where
$(\cdot , \cdot )_{h}$ denotes that this term is evaluated using an appropriate
numerical integration rule
to account for the superconvergence of $\nabla p_h$, and
define the following local post-processing technique

\paragraph{Problem LPP} {\it For a given $p_h$, solution of { Problem}
{Ph}, find $\bm{u}_h \in \mathbb{M}_{h}$ such that
$$
 b_{h} (\bm{u}_h, \bm{v}_h ) = ( \nabla p_h, \bm{v}_h)_{h} -
 h^2 (f, \Div\bm{v}_h)  \quad \forall \bm{v}_h \in\mathbb{M}_{h}  \, .
$$
}
Stability and convergence of this local post-processing formulation is
proved in \cite{LOULA95}. For quadrilateral elements and $k>1$, optimal orders of
convergence are proved for any macroelement configuration, including element by element
post-processings, when the mesh-dependent term $(\cdot,\cdot)_h$
is computed using the $k\times k$ Gauss integration points where the
finite element approximation of the gradient is superconvegent with
$$
  | \nabla p - \nabla p_h |_h   \leq  C  h^{k+1} \| p \|_{k+2} \, ,
  \quad  | \nabla q |_h = (\nabla q , \nabla q)^{1/2}_{h} \, .
$$
For $k=1$ this
post-processing is also stable and optimally convergent for any macroelement composed
by at least two adjacent elements with a common edge. In these cases the bilinear
form $ b_{h} (\cdot , \cdot)$ defines the norm
$$
||| \bm{u}_h ||| = ( b_{h} (\bm{u}_h , \bm{u}_h) )^{1/2}
       \quad \forall \bm{u}_h \in \mathbb{M}_{h} \, ,
$$
which is equivalent to the $[L^2(\Omega)]^2$-norm on $\mathbb{M}_{h}$.
By usual arguments, we derive the estimate
$$
||| \bm{u} - \bm{u}_h ||| \le C  ( ||| \bm{u} - \bm{v}_h  ||| +
    | \nabla p - \nabla p_h |_h )   \quad \forall \bm{v}_h \in \mathbb{M}_{h} \, .
$$
For regular solutions,
the above estimate leads to the optimal orders of convergence
\begin{equation}
\label{eq:macro-l2}
 \| \bm{u} - \bm{u}_h \| \leq  C  h^{k+1} \| \bm{u} \|_{k+1} \, ,
\end{equation}
for the post-processed velocity approximation, in $[L^2(\Omega)]^2$-norm.
For  problems with regular exact solutions  concerning
homogeneous and isotropic  media, these optimal orders of convergence
are confirmed in \cite{LOULA95} in a large number of numerical experiments.
An extension of this local post-processing for triangular elements
is presented in \cite{RITA2002}, in the context of h-adaptive
 analysis of Poisson's problem,
using the superconvergence points  identified by
Babuska, Strouboulis, Upadhyay and Gangaraj  in \cite{BABUSKA96}.

\begin{figure}[H]
\centering
\subfloat[]{\includegraphics[angle=90,height=.3\textwidth,angle=0]{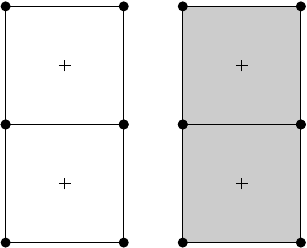}
\label{fig:sm1}}
\hspace{1cm}
\subfloat[]{\includegraphics[angle=90,height=.3\textwidth,angle=0]{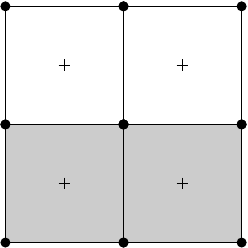}
\label{fig:sm2}}
\caption{Stable macroelement configurations composed of bilinear
    quadrilaterals ($\bullet$, nodal points; +, superconvergent points):
    (a) two homogeneous macroelements composed by an arrange of
    $2\times1$ bilinear quadrilaterals; (b) a macroelement of
    $2\times2$ elements with interface of discontinuity. In this last
    case, the method described in Section \ref{sec:ppge} is adopted
    inside the macroelement} \label{fig:macros}
\end{figure}

Our point here is the application of this local post-processing to
an anisotropic porous media with an interface of material
discontinuity.  Contrary to the global $C^0$ post-processing, the
local post-processing is posed in finite element subspaces of
$[L^2(\Omega)]^2$ using Lagrangian interpolation which are
discontinuous on the macroelement interfaces. Applied to porous
media with an interface of discontinuity on the hydraulic
conductivity, the local post-processing naturally  accommodates
discontinuities of the tangent component but it does  not ensure
continuity of the normal component of the velocity field on the
macroelement interfaces. However, as it is a stable and
variationally consistent formulation, we should expect optimal orders
of convergence for sufficiently regular exact solutions on each
subdomain $\Omega_i$, independently of the fact that the
continuity/discontinuity constraints on the normal/tangential
components of Darcy's velocity are not exactly imposed on the
interface of two media with different hydraulic conductivity.
Numerical results presented in Sec. \ref{sec:ex4} confirm the optimal
orders of convergence when the interface of discontinuity belongs to
the interface of the macroelements.  If  the interface of
discontinuity is put in the interior of a macroelement, the
continuity/discontinuity constraints on the normal/tangential
components of Darcy's velocity can be exactly imposed using the
method presented in Sec. \ref{sec:ppge}. In Figure \ref{fig:macros}
we illustrate: (a) case where the interface belongs to the edges of two
homogeneous macroelements, (b)  case where the interface is put in
the interior of a macroelement. Again, optimal orders of convergence
should be expected for piecewise regular solutions. This fact is
also confirmed numerically in Sec. \ref{sec:ex4}.

%===================================================================
%
%
\section{Numerical Examples}
\label{sec:ex}
%
%
%===================================================================
%

In this section, we present  examples of  applications of the proposed
post-proces\-sing techniques to Darcy flow in layered heterogeneous
media and compare their results with stabilized mixed methods. The
first three examples illustrate qualitative aspects of the
approximations, while the last example test the predicted orders of
convergence. In all simulations we adopted for the global post-processings,
with or without interface of discontinuity,
$\delta=1$ and $\alpha=1$. According to the numerical analysis
presented in \cite{LOULA95} $\alpha =1$ is the best choice leading
to quasi optimal orders of convergence for the velocity approximation
in $L^2$-norm, and $\delta>0$ is sufficient for stability.

\subsection{Flow Between Parallel Plates}
\label{sec:ex1}
\begin{figure}[h]
\centering
\includegraphics[width=.6\linewidth]{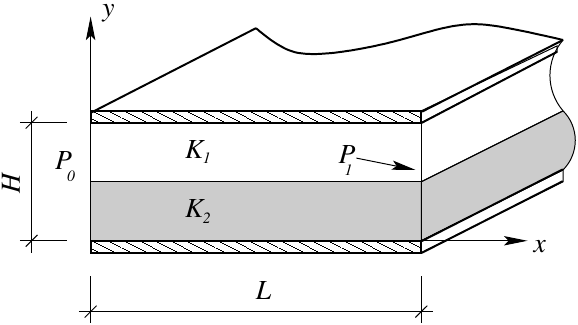}
  \caption{Domain.}
  \label{fig:pp-fig}%
\end{figure}

This simple example clearly shows the consequences of imposing $C^0$
continuity in problems with discontinuous solutions. Let
$\Omega=[0,L]\times[0,H]$ be the longitudinal section of a
heterogeneous porous medium composed by two regions with different
conductivities, confined between two parallel plates
of infinite extent, as shown in Fig \ref{fig:pp-fig}.% 4.
\begin{figure}[h]
\centering
  \includegraphics[height=.5\linewidth]{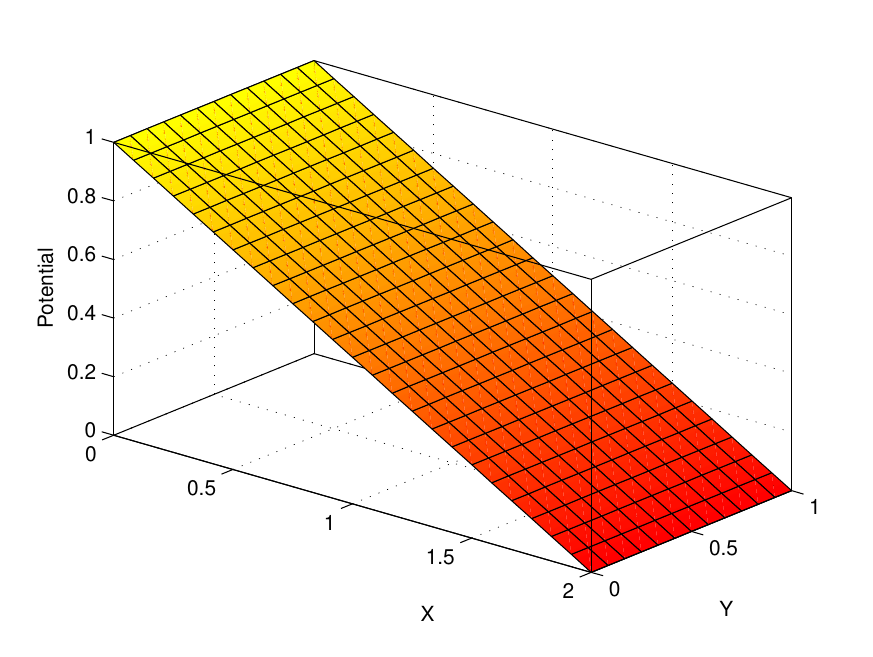}
  \caption{Approximated potential.}
  \label{fig:pp-fig-con}%
\end{figure}
%\pagebreak
%
%%%
On the left extremity the hydraulic potential is prescribed
$P_0=1.0$m and on the right extremity $P_1=0.0$m; on the plates we
consider no-flux boundary conditions $\bm{u}\cdot\bm{n}=0$. The dimensions
adopted are $L=2.0$m and $H=1.0$m, and the conductivities are
$K_1=2.0\mathrm{m/day}$ and $K_2=1.0\mathrm{m/day}$.
The flow is horizontal, originated by the gradient of
potential $dp/dx=-0.5$, and the velocity can easily be calculated
using Darcy's law (\ref{eq:darcy}), yielding  piecewise constant
velocities $u_1=1.0\mathrm{m/s}$ in the upper region and
$u_2=0.5\mathrm{m/s}$ in the lower one.
%
%%%

\begin{figure}[htb]
\centering
\subfloat[]{\includegraphics[angle=0,width=.49\textwidth,angle=0]{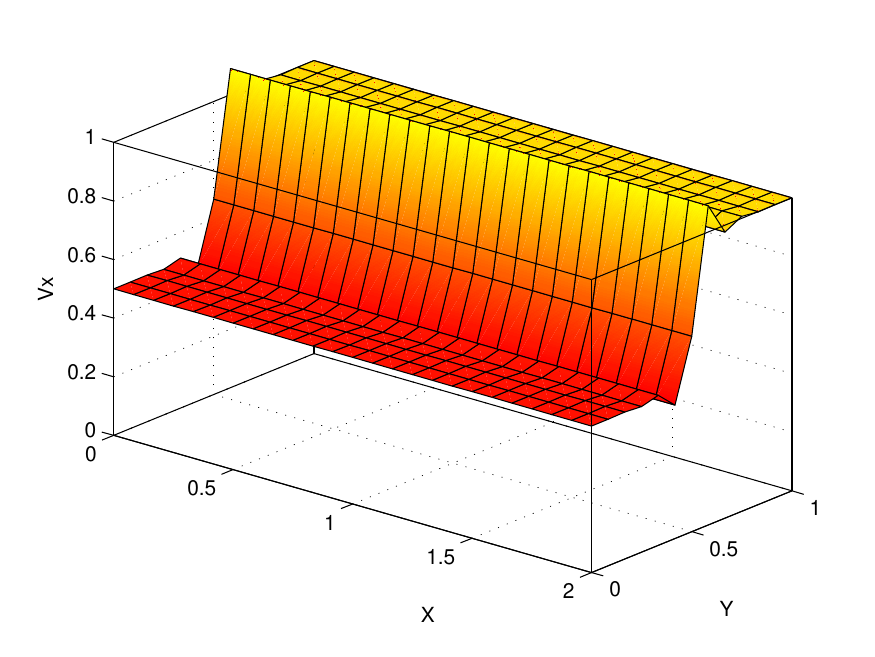}
\label{fig:pp1}}
%
% \hspace{1cm}
%
\subfloat[]{\includegraphics[angle=0,width=.49\textwidth,angle=0]{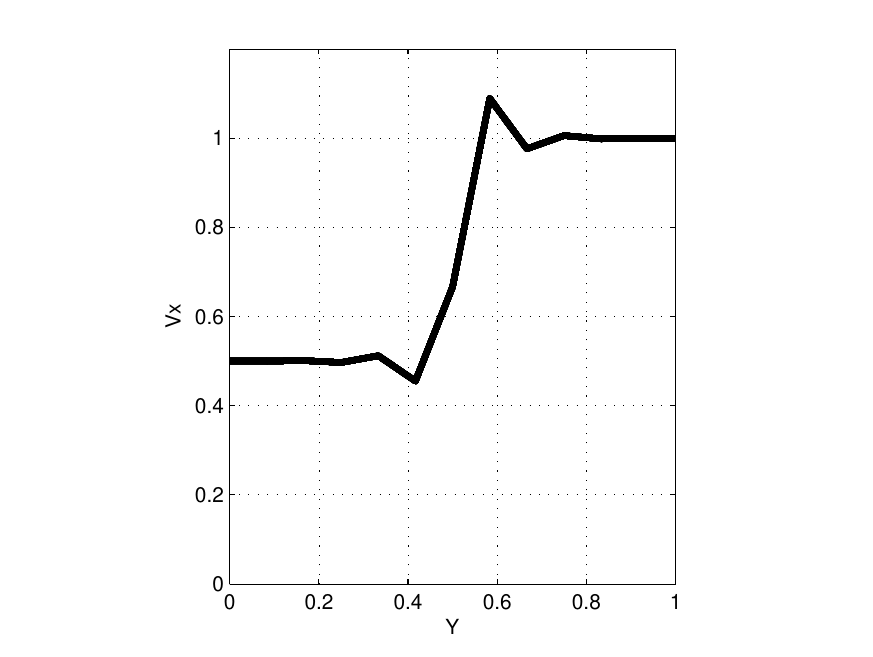}
\label{fig:pp2}}
\caption{(a) HVM, GLS and GPP formulations: component $u_x$ of the approximated
    velocity. (b) Projection in $yz$ plane.}
    \label{fig:pp-sem}%
\end{figure}

\begin{figure}[htb]
\centering
\subfloat[]{\includegraphics[angle=0,width=.49\textwidth,angle=0]{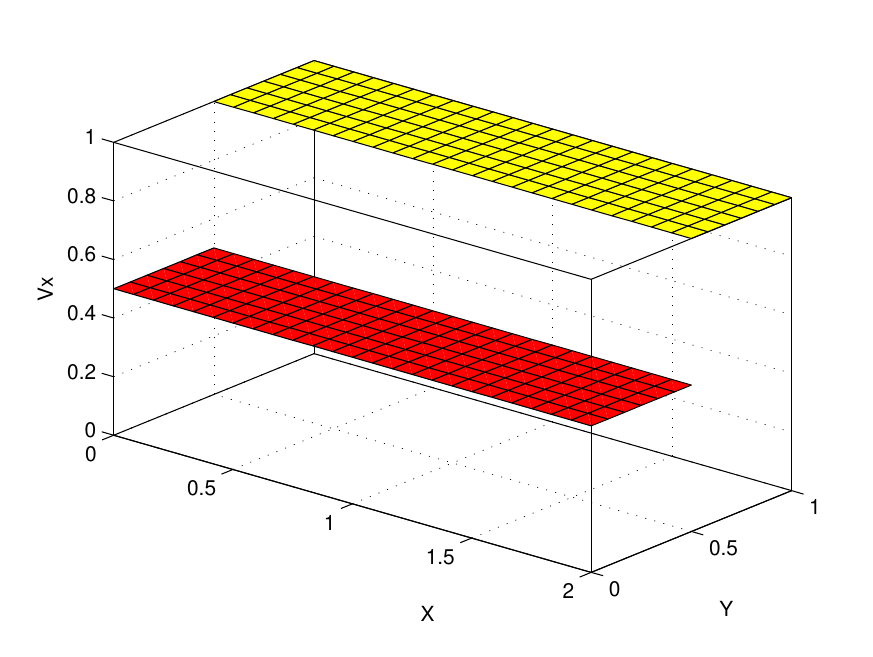}
\label{fig:pp3}}
%
% \hspace{1cm}
%
\subfloat[]{\includegraphics[angle=0,width=.49\textwidth,angle=0]{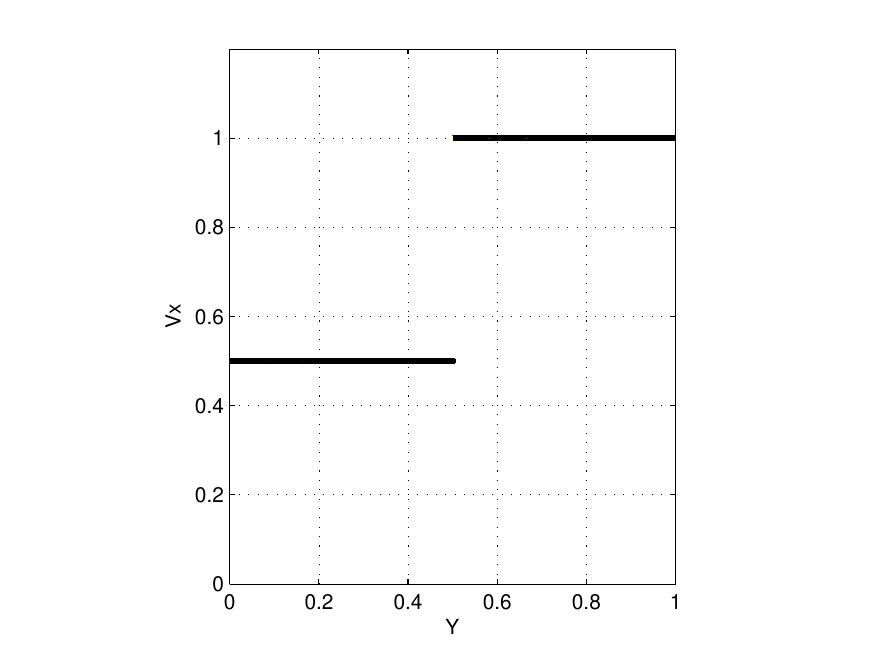}
\label{fig:pp4}}
\caption{(a) Velocity approximated with GPPID, LPP and LPP with incorporated interface.
(b) Projection in $yz$ plane.}
    \label{fig:pp-com}%
\end{figure}

First we consider $C^0$ approximations by the stabilized mixed methods:
HVM and GLS (with $\delta_1=\delta_2=1/2$), and by the global $C^0$ post-processing
(GPP). A uniform mesh of $24\times 12$ bilinear quadrilateral elements is
adopted in all experiments.
No difference was observed in the solutions obtained with these methods
using $C^0$ velocity interpolations.
Approximations for potential and velocity
fields are shown in Figures \ref{fig:pp-fig-con} and
\ref{fig:pp-sem}, respectively. We can clearly observe spurious
oscillations in the velocity approximation, as a consequence of
imposing continuity. In contrast, in Fig \ref{fig:pp-com}
we present continuous/discontinuous  velocity approximations obtained with:
the global post-processing incorporating the discontinuity (GPPID),
the local post-processing with $12\times 6$ homogeneous macroelements
composed by arranges of $2\times 2$ bilinear quadrilaterals
and the local post-processing with $6\times 3$ macroelements composed by arranges
of $4\times 4$ elements, in which discontinuity/continuity constraints
are introduced at macroelement level.
The solution now is in perfect agreement with the real flow.
\subsection{Flow With Less Pervious Barriers}
\label{sec:ex2}

\begin{figure}[htb]
  \centering
  \includegraphics[width=.60\linewidth]{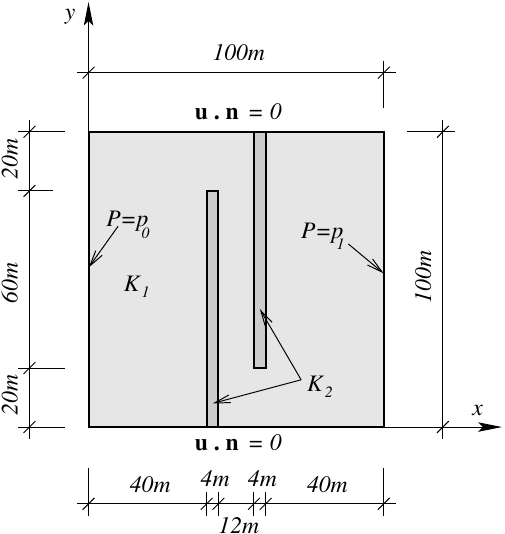}
  \caption{Domain.}
  \label{fig:tese-dominio}
\end{figure}

This example, from Mos\'e et al. \cite{MOSE94}, simulates the flow
in a porous media intersected by two less pervious barriers. These
barriers force the flow to pass through a kind of channel. The
boundary conditions are the same of the previous example. The
geometry is detailed in  Fig \ref{fig:tese-dominio}. The
conductivities are $K_1=1.0$m/day for the porous media and
$K_2=1.0\times10^{-5}$m/day for the barriers.

\begin{figure}[htb]
  \centering
  \includegraphics[width=.55\linewidth]{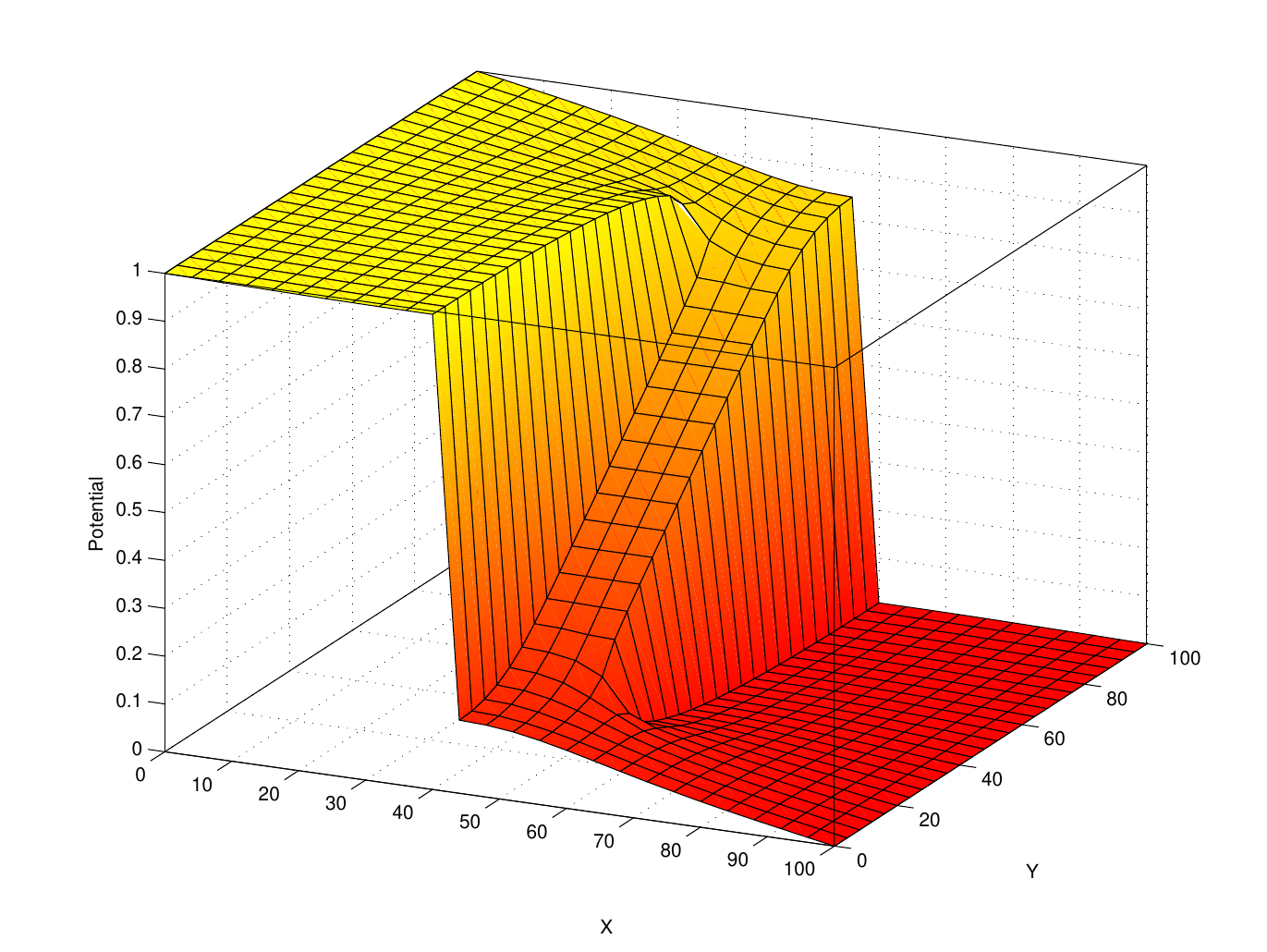}
  \caption{Approximated potential.}
  \label{fig:tese-fig-con}
\end{figure}

We approximated this problem by HVM and GPPID methods, using  a
uniform mesh of $25\times 25$ bilinear quadrilaterals. For GPPID, we
assumed that at the corners with material discontinuity the relation
$K_1^{-1}\bm{u}_1=K_2^{-1}\bm{u}_2$ holds. This is certainly a crude
approximation for the interface conditions in which only the
restriction on the tangential component is imposed. We should remind
that in our model problem the interface is assumed to be smooth.
This smoothness hypothesis is violated in this example.
\begin{figure}[htb]
\centering
\subfloat[HVM.]{\includegraphics[angle=0,width=.49\textwidth,angle=0]{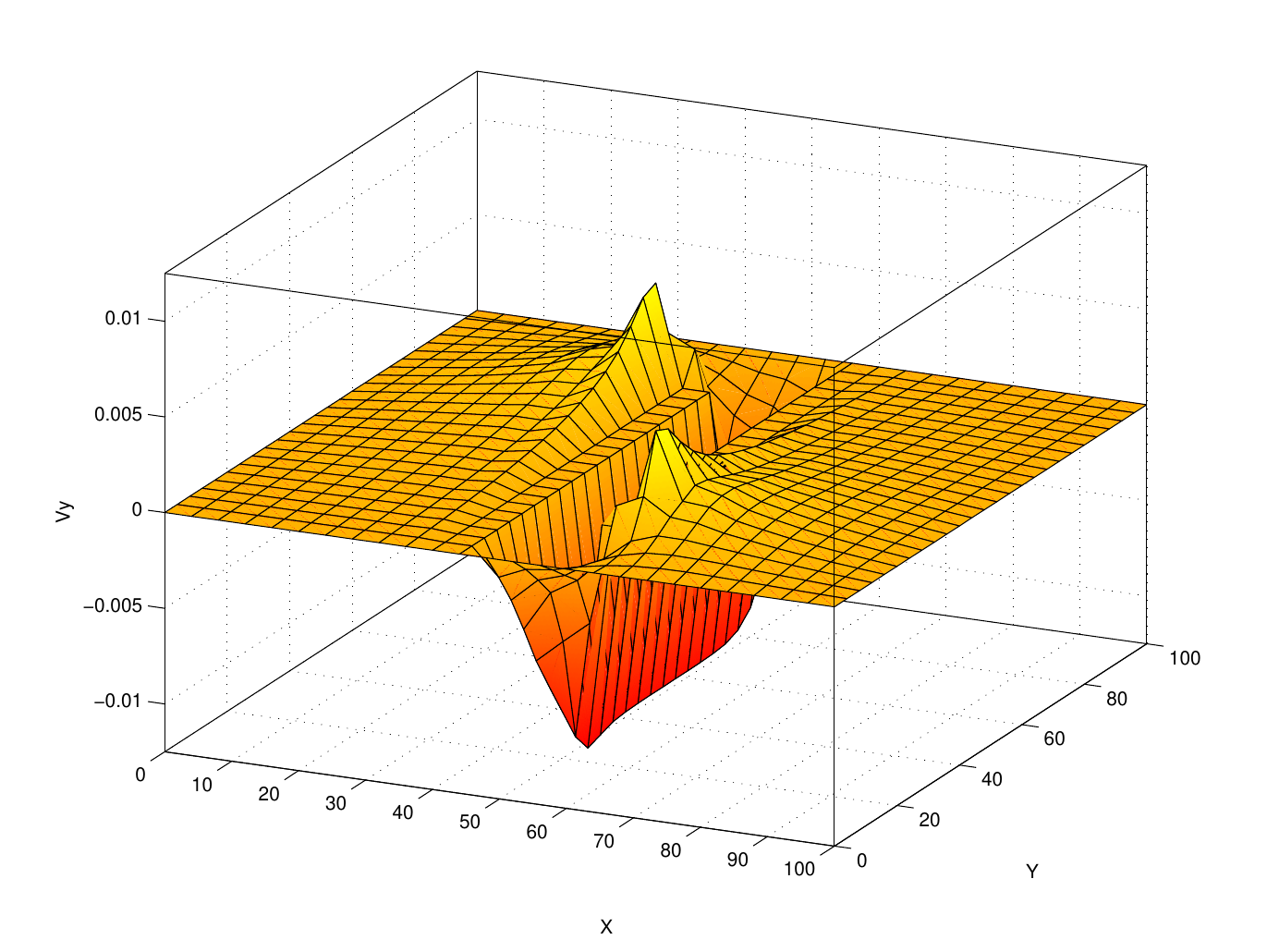}
\label{fig:bidu-vy-hvm}}
%
% \hspace{1cm}
%
\subfloat[GPPID.]{\includegraphics[angle=0,width=.49\textwidth,angle=0]{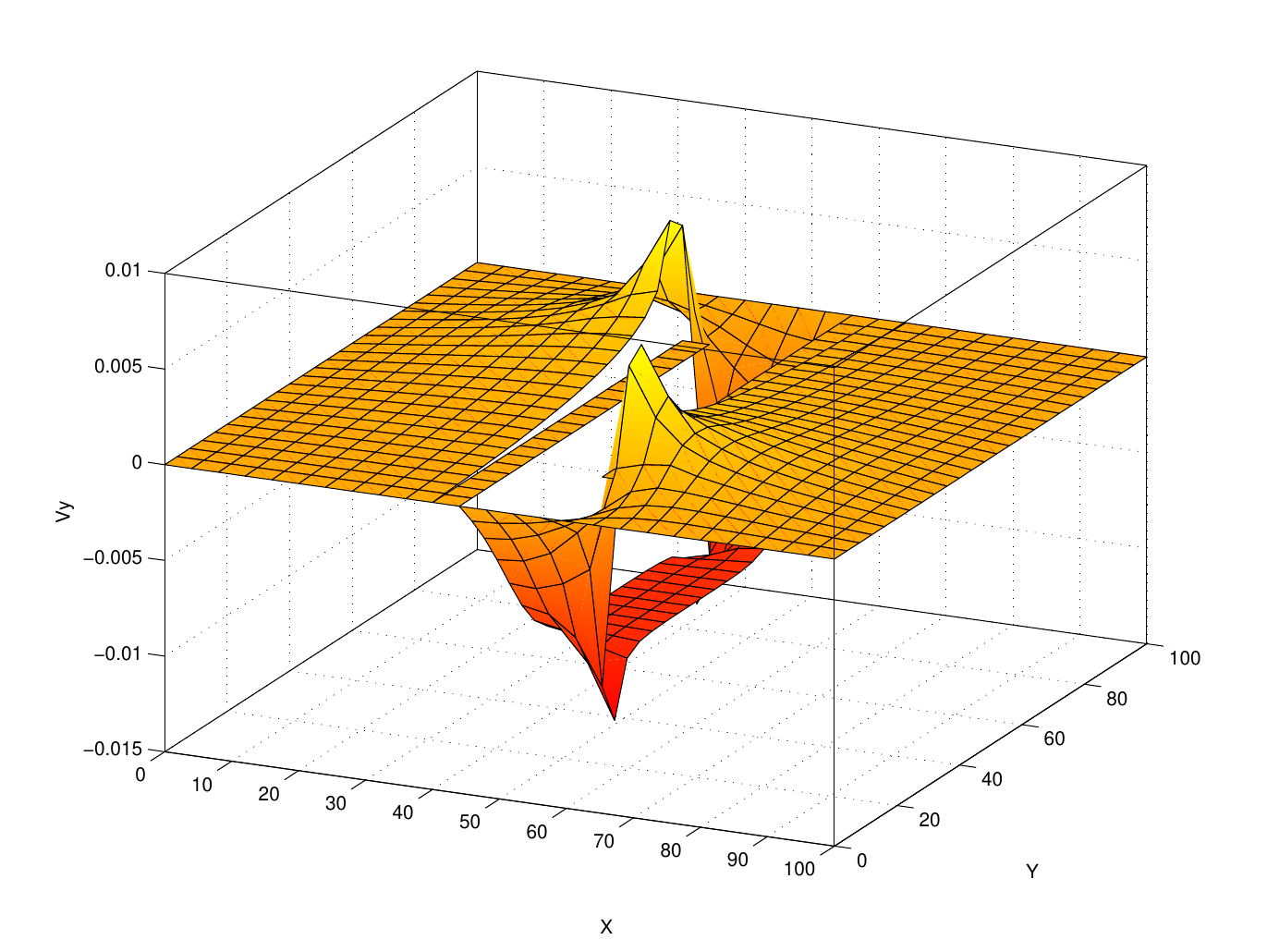}
 \label{fig:bidu-vy-gppid}}
\caption{Flow with less pervious barriers. Component $u_y$.}
    \label{fig:bar-com}%
\end{figure}
\begin{figure}[htb]
\centering
\subfloat[Flow by HVM (There is artificial adherence on interface).]{\includegraphics[angle=0,width=.49\textwidth,angle=0]{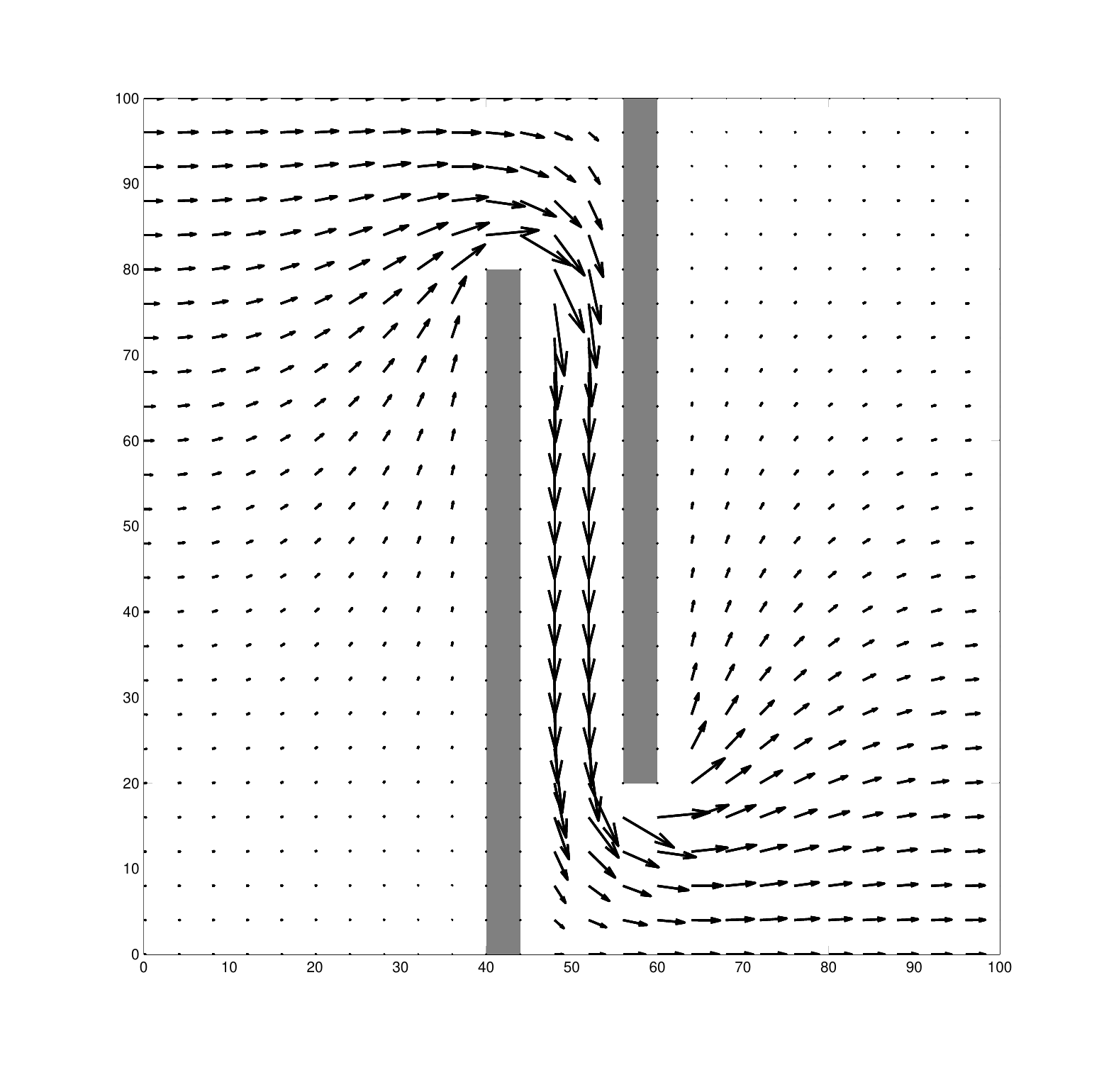}
\label{fig:bidu-v-hvm}}
%
% \hspace{1cm}
%
\subfloat[Flow by GPPID (Tangent component well represented
      on interface).]{\includegraphics[angle=0,width=.49\textwidth,angle=0]{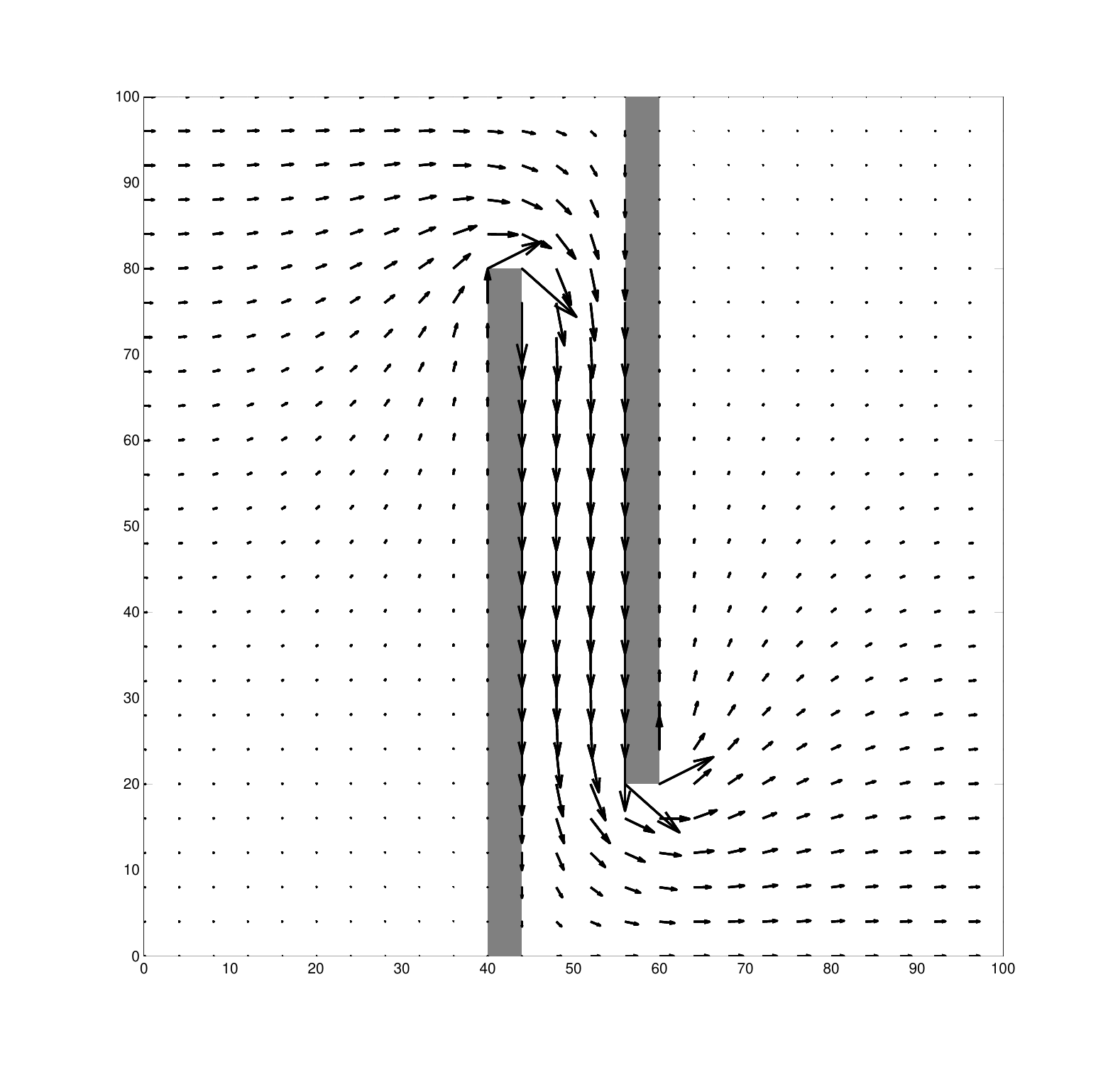}
 \label{fig:bidu-v-gppid}}
\caption{Flow with less pervious barriers. Flow fields.}
    \label{fig:bar-com2}%
\end{figure}
Figure \ref{fig:tese-fig-con} shows the Galerkin approximation for
the potential. The potential approximated with HVM is similar.
Figures \ref{fig:bidu-vy-hvm} and \ref{fig:bidu-vy-gppid} show the
$u_y$ velocity component obtained with HVM and GPPID,  respectively.
We observe that the GPPID  approximation capture the discontinuity.
The flow calculated with HVM exhibits an artificial adherence to the
barriers while the GPPID solution presents a  behavior typical of
no-flux condition, as we clearly observe in Figures
\ref{fig:bidu-v-hvm} and \ref{fig:bidu-v-gppid}.

\subsection{Porous Media with Three Different Conductivities}
\label{sec:ex3}
This is another example that does not fit in our model problem.
Again we do not have a smooth interface of discontinuity.
It is concerned with a heterogeneous porous media composed by three
subdomains with different conductivities and two intersecting
interfaces, as suggested by one of the referees. At the point of
intersection it is not possible to assemble the nodal transformation
(\ref{eq:uref2}). To handle this problem  we adopted the local
post-processing formulation using macroelements, naturally
accommodating the interfaces of discontinuity on the macroelement
edges or imposing the continuity/discontinuty constraints, through
the nodal transformation (\ref{eq:uref2}), on interfaces interior to
the macroelements.

\begin{figure}[htb]
\centering
\subfloat[]{\includegraphics[angle=0,width=.32\textwidth,angle=0]{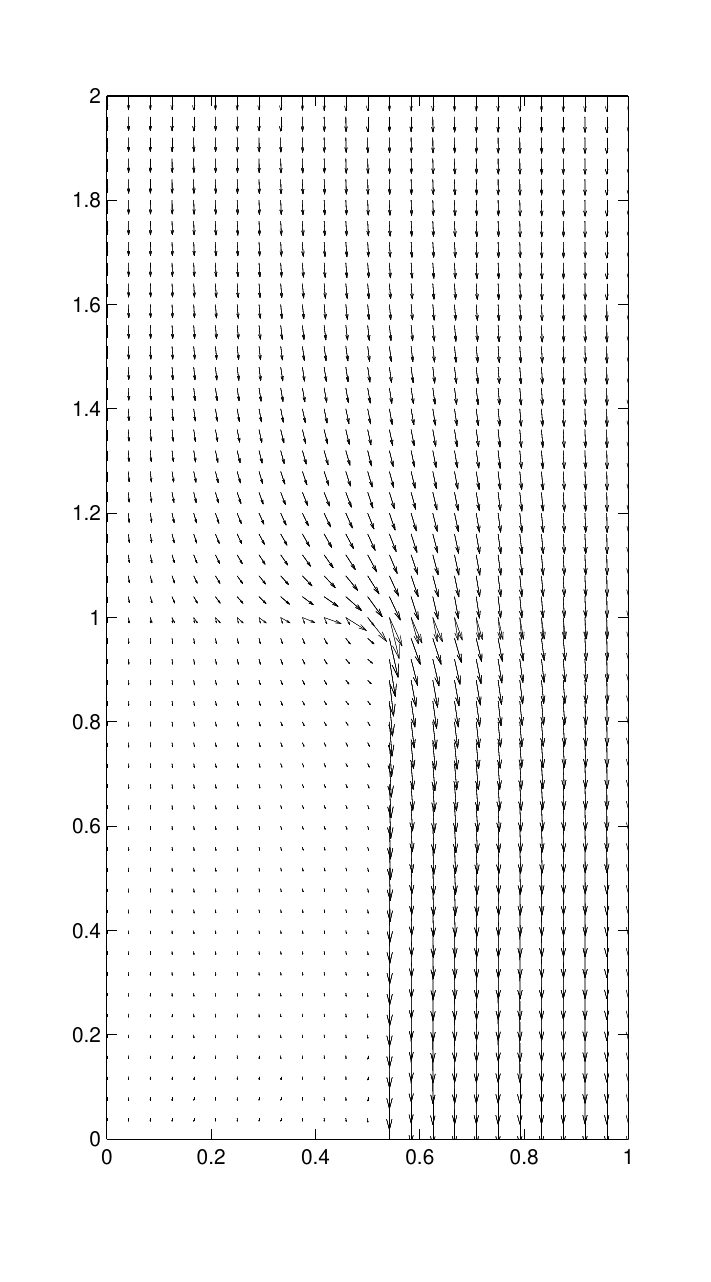}
\label{fig:macrosa}}
%
% \hspace{1cm}
%
\subfloat[]{\includegraphics[angle=0,width=.32\textwidth,angle=0]{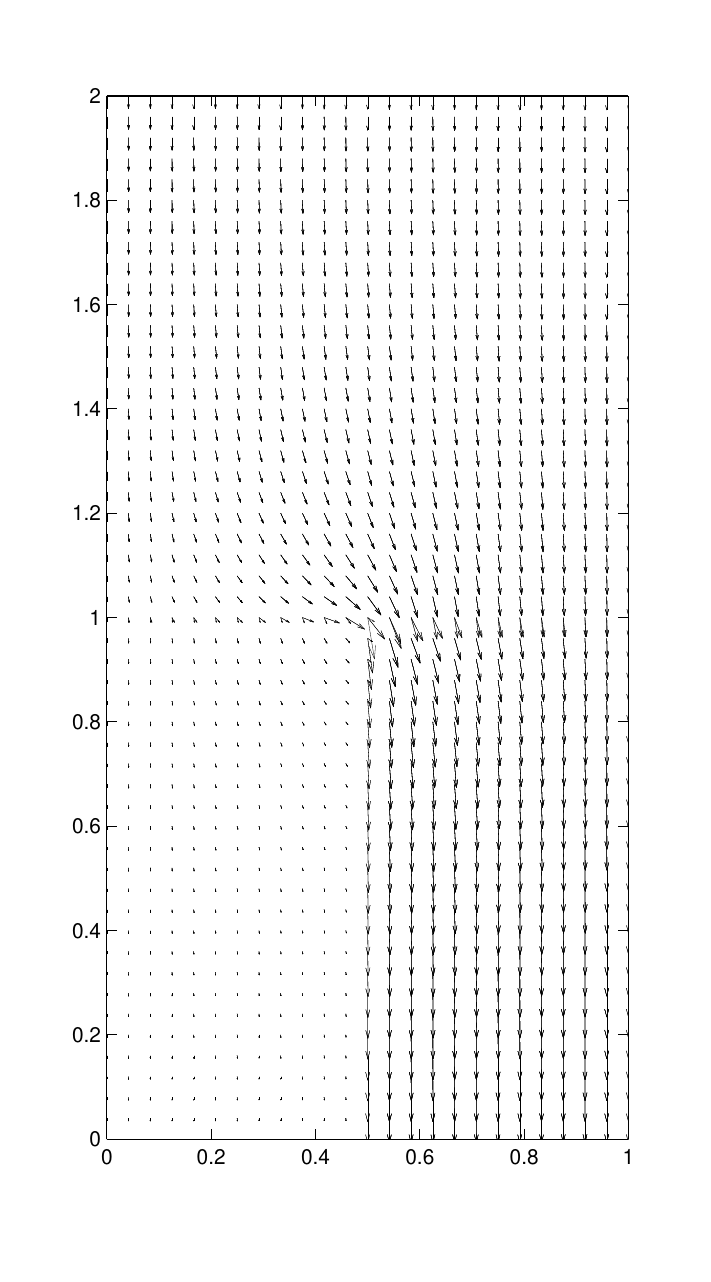}
    \label{fig:macrosb}}
%
% \hspace{1cm}
%
\subfloat[]{\includegraphics[angle=0,width=.32\textwidth,angle=0]{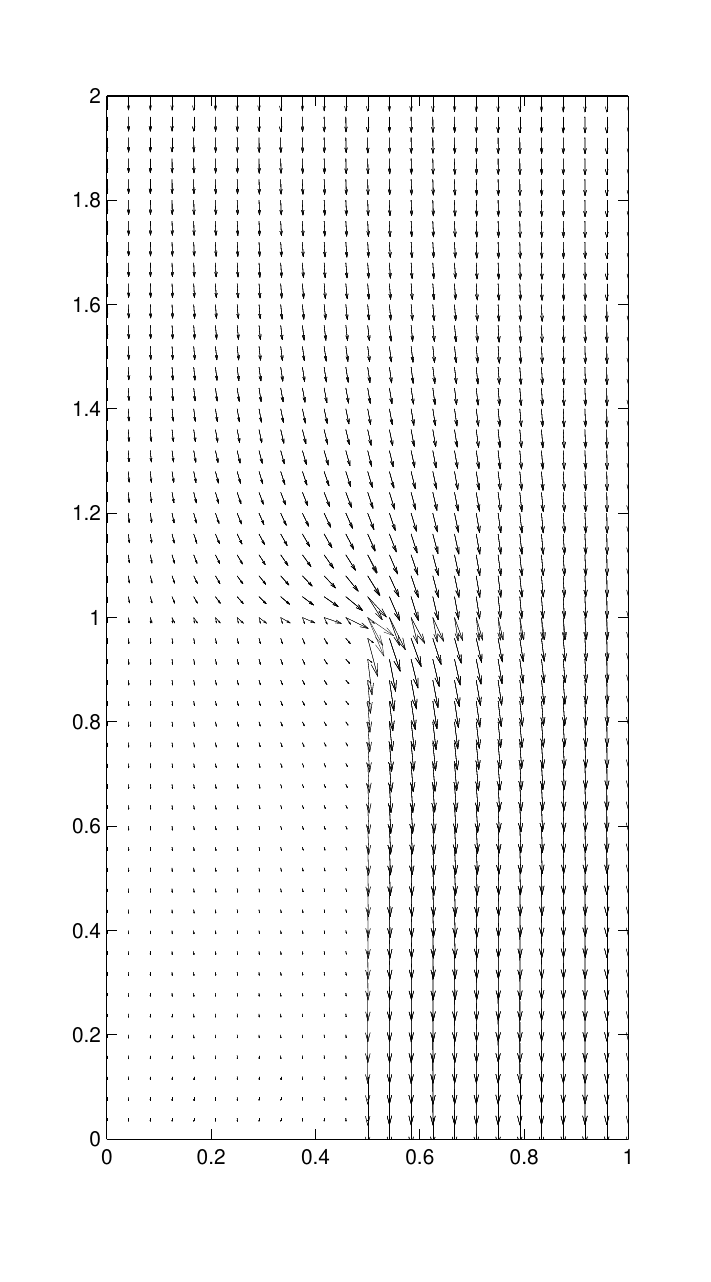}
    \label{fig:macrosc}}
\caption{Flow through the porous media with three different
conductivities. (a) Two (continuous)
macroelements of $24 \times 25$ elements, (b) interface in the interior of the first macroelement
incorporated and (c) four macroelements of $12\times 25$ elements.}
    \label{fig:macrosf}%
\end{figure}
The domain of the problem is a rectangle given by the intervals
$[0.0,1.0]\times[0.0,2.0]$ (lengths in meters), composed by three
subdomains: the first defined on $[0.0,0.5]\times[0.0,1.0]$ with
conductivity $K=1\mathrm{m/day}$, the second on
$[0.5,1.0]\times[0.0,1.0]$ with $K=10\mathrm{m/day}$ and the third
on $[0.0,1.0]\times[1.0,2.0]$ with $K=5\mathrm{m/day}$. The
hydraulic potential is prescribed $P_1=1.0\mathrm{m}$  on the upper
extremity and $P_0=0.0\mathrm{m}$ at the bottom; at the vertical
boundaries we consider no-flux boundary conditions $\bm{u}\cdot\bm{n}=0$.
The approximations were performed with a uniform mesh of
$24\times50$ bilinear quadrilaterals, composing the following
arranges: (a) two macroelements of $24\times 25$ elements, where the
interface between the first and the second subdomains is interior to
the first macroelement, (b) four macroelements of $12\times25$
elements; in this last case all interfaces belong to the set of
edges of the macroelements.

\begin{figure}[htb]
\centering
\subfloat[]{\includegraphics[angle=0,width=.32\textwidth,angle=0]{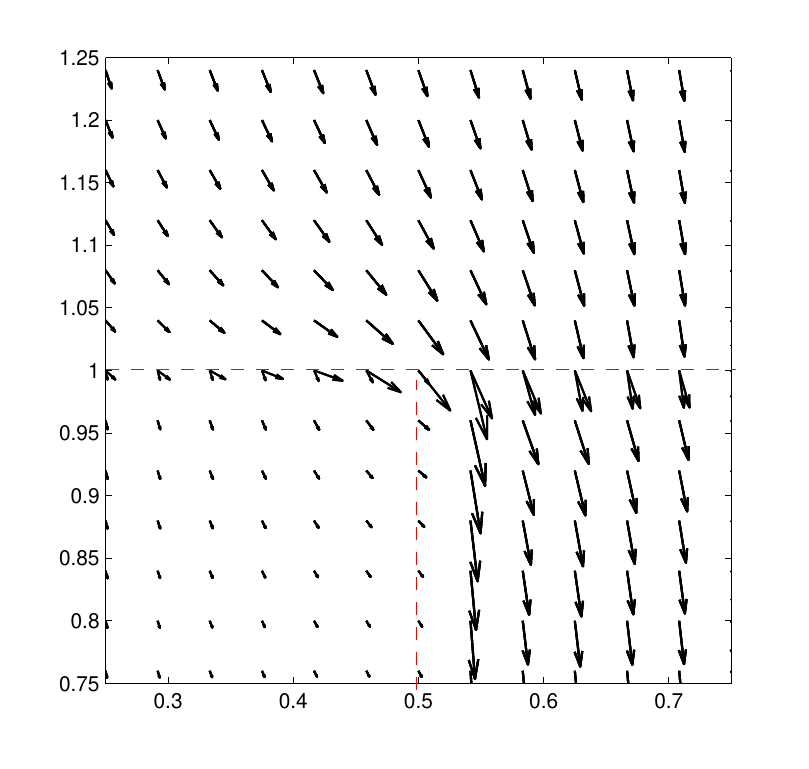}
\label{fig:macrosza}}
%
% \hspace{1cm}
%
\subfloat[]{\includegraphics[angle=0,width=.32\textwidth,angle=0]{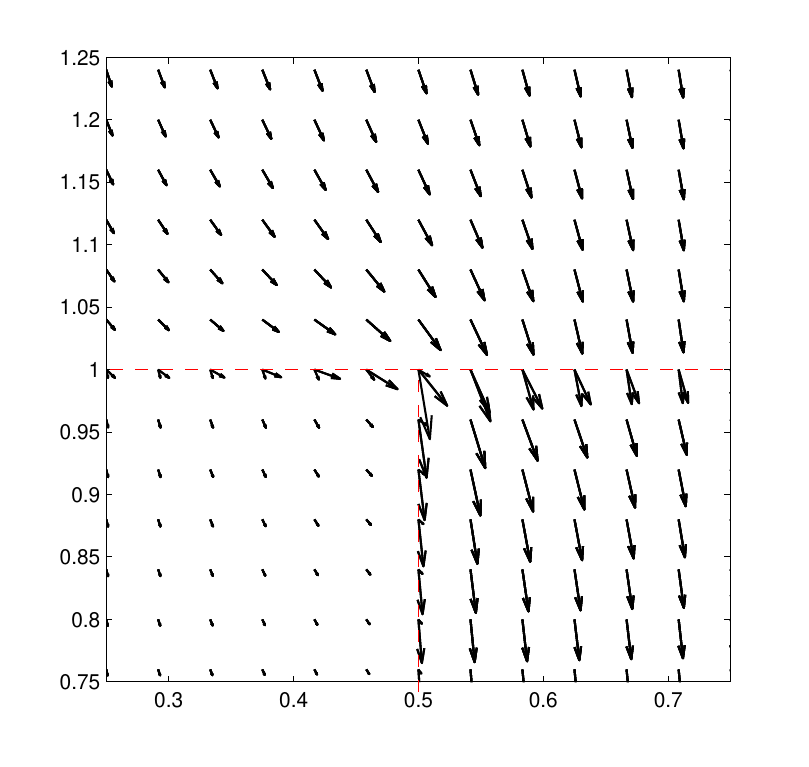}
\label{fig:macroszb}}
%
% \hspace{1cm}
%
\subfloat[]{\includegraphics[angle=0,width=.32\textwidth,angle=0]{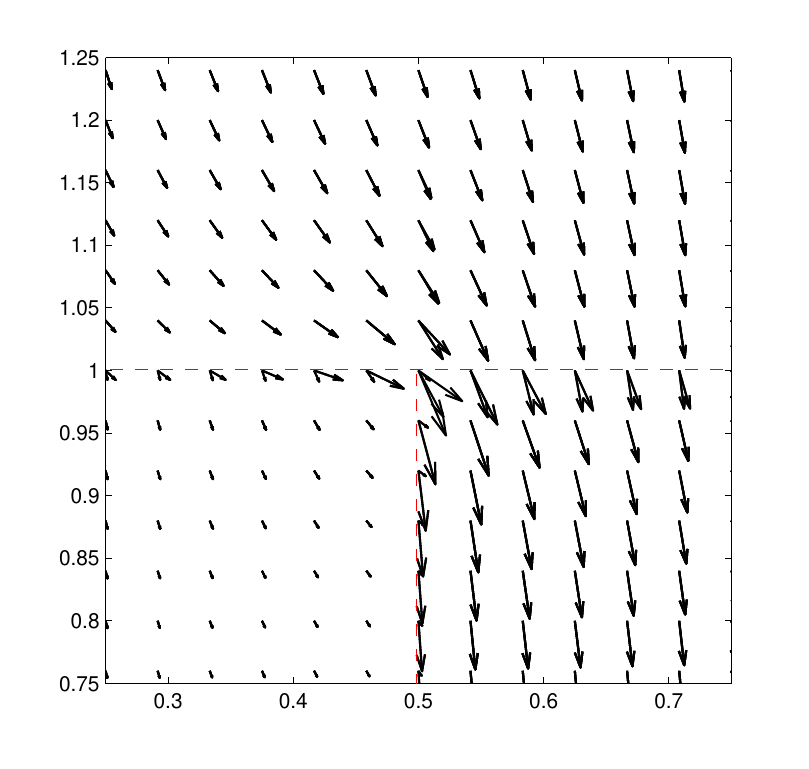}
\label{fig:macroszc}}
\caption{Details of the flow near to the point of intersection
    of the interfaces. (a) Two (continuous)
      macroelements,
      (b) two macroelements (ID) and (c) four macroelements.}
    \label{fig:macrosz}%
\end{figure}

In Figure \ref{fig:macrosa} we show the flow for the approximation
with 2 macroelements. In this case, the solution is continuous
inside each macroelement.  Adherence of the flow on the interface
between the lower subdomains can be observed. A significant
improvement on the flow approximation is achieved when we impose the
interface conditions as shown in Fig \ref{fig:macrosb}. The velocity
on the interface is recovered. Similar results are obtained with the
arrange of 4 macroelements with natural interface condition as
presented in Fig \ref{fig:macrosc}. Details of the flow near to the
point of intersection of the interfaces are shown in Figures
\ref{fig:macrosza}, \ref{fig:macroszb}, and \ref{fig:macroszc}.
\subsection{Convergence Study: An Anisotropic Heterogeneous Example}
\label{sec:ex4}

\begin{figure}[htb]
\centering
\subfloat[HVM.]{\includegraphics[angle=0,width=.49\textwidth,angle=0]{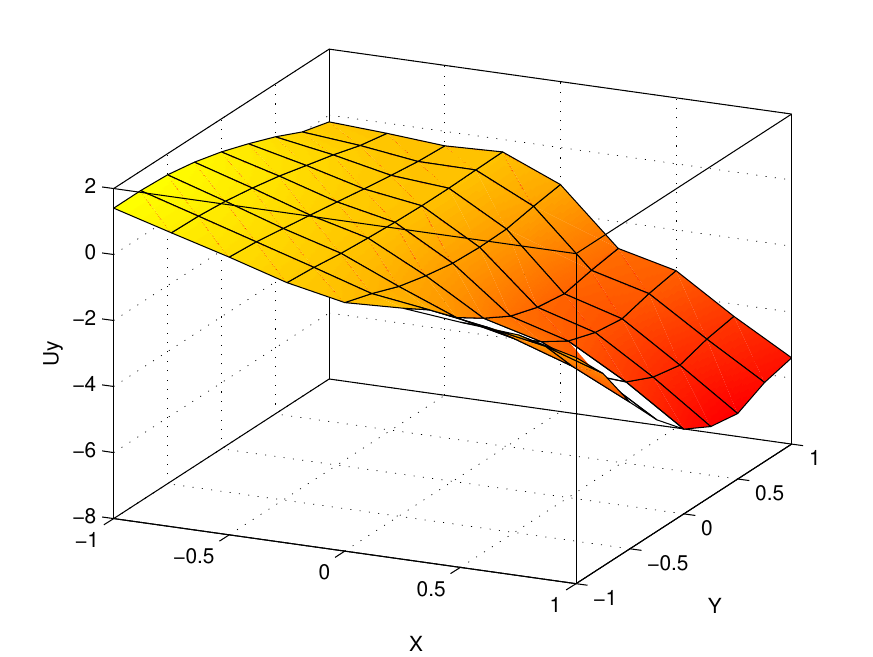}
\label{fig:vy-anis-hvm}}
%
% \hspace{1cm}
%
\subfloat[GPPID.]{\includegraphics[angle=0,width=.49\textwidth,angle=0]{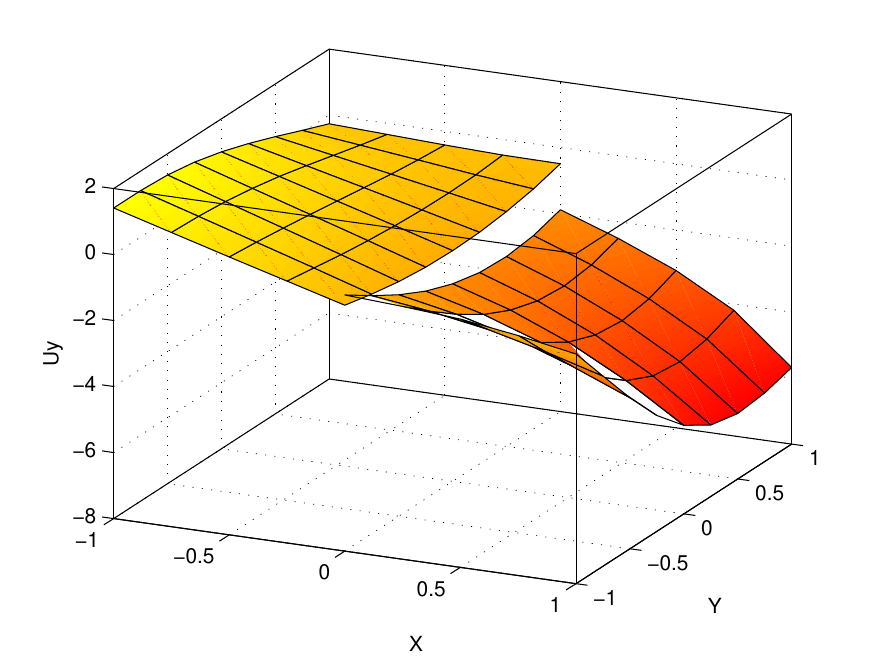}
 \label{fig:vy-anis-gppid}}
\caption{Convergence study. Component $u_y$.}
    \label{fig:bar-conv}%
\end{figure}

In this example we perform  a convergence study in a general case
where the conductivity tensor is anisotropic and heterogeneous. This
test problem from Crumpton, Shaw and Ware \cite{CRUMPTON95} is
defined on the square $[-1,1]\times[-1,1]$, with Dirichlet boundary
conditions. The conductivity is given by
  \[
  K= \left(
      \begin{array}{ccc}
    1 &  0  \\
    0 &  1
      \end{array}
      \right)   \quad  x<0,  \qquad
   K=\gamma\left(
      \begin{array}{ccc}
    2 &  1  \\
    1 &  2
      \end{array}
      \right)   \quad x>0,
  \]
where the parameter $\gamma$ is used to vary the strength of the
discontinuity at $x=0$. In our experiments we used $\gamma=1.0$.
The exact potential field is given by
  \[
  p=
  \left\{
  \begin{array}{lcr}
  \gamma (2\sin(y)+\cos(y)) x + \sin(y), & & x<0 \\ \\
   \exp(x)\sin(y),                       & & x>0
  \end{array}
  \right.
  \]
and the velocity field can be directly calculated by Darcy's law. As
in all formulations considered here the potential belongs to $\mathcal{Q}\in
H^1(\Omega)$, we prescribe this variable on the boundary. The finite
element solutions were computed adopting uniform meshes of
$8\times8$, $16\times16$, $32\times32$ and  $64\times64$ bilinear
quadrilaterals (Q1) and $4\times4$, $8\times8$, $16\times16$ and
$32\times32$ biquadratic quadrilaterals (Q2). In the specific
convergence study for the local post-processing with interface of
discontinuity in the interior of macroelements, we adopted uniform
meshes of $6\times 6$, $12\times 12$, $24\times 24$ and $48\times
48$ Q1 elements,  and $6\times 6$, $12\times 12$ and $24\times 24$
Q2 elements. To illustrate the different behaviors of continuous and
discontinuous approximations, the $y$ components of the velocity
fields obtained with $8 \times 8$ bilinear quadrilaterals using HVM
and GPPID are shown in Figures \ref{fig:vy-anis-hvm} and
\ref{fig:vy-anis-gppid}, respectively.

According to (\ref{eq:est2}), for sufficiently
regular solution we should expect the following orders of
convergence
\[
\| \bm{u}- \bm{u}_h \| \leq Ch^{k+0.5}, \quad \| \Div \bm{u}- \Div \bm{u}_h \|
\leq Ch^k,
\]
for the global $C^0$ post-processing with $\alpha=1$ and $\delta>0$.
The error estimates for the other methods are given by Equations
(\ref{eq:mpe-p}-\ref{eq:mpe-u}) for GLS, (\ref{eq:hvm-p}-\ref{eq:hvm-u}) for HVM and
(\ref{eq:macro-l2}) for LPP. The point here is that the exact solution of this problem
is not sufficiently regular in the whole domain $\Omega$, but  regular
enough in each subdomain $\Omega_1$ and $\Omega_2$, so that we can expect
confirmation of the predicted orders of convergence  whenever the discontinuity
in the interface is appropriately captured.
In the convergence graphics, we present plots of $L^2(\Omega)$ norm
of the errors of the potential, velocity and its divergence versus
$-\log(h)$.
The label MACRO 2X2 (ID) in the pictures identifies approximations
obtained with local post-processing with macroelements composed by
$2\times2$ elements with interface of discontinuity in the interior
os macroelements (exactly imposed constraints),
while MACRO 2X2 corresponds to approximations obtained with
macroelements composed by $2 \times 2$ homogeneous elements with the
interfaces of discontinuties belonging to to the set of macroelement
edges (constraints naturally accommodated). Again, we have adopted
$\delta_1=\delta_2=1/2$ in the GLS formulation.

\begin{figure}[htb]
\centering
\subfloat[Bilinear.]{\includegraphics[angle=0,width=.49\textwidth,angle=0]{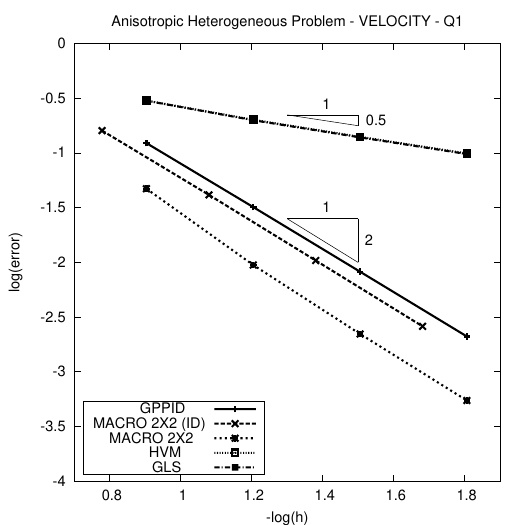}
\label{fig:vs-q1}}
%
% \hspace{1cm}
%
\subfloat[Biquadratic.]{\includegraphics[angle=0,width=.49\textwidth,angle=0]{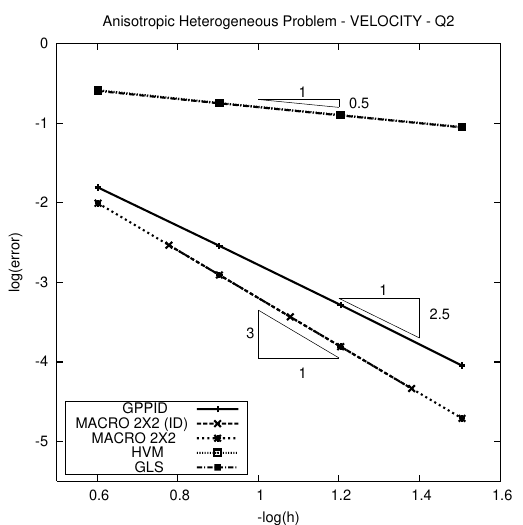}
 \label{fig:vs-q2}}
\caption{Convergence of the velocity field.}
    \label{fig:bar-conv-v}%
\end{figure}

Convergence results of velocity approximations for Q1 and Q2
elements are presented in Figures \ref{fig:vs-q1} and \ref{fig:vs-q2},
respectively.
Convergence of $O(h^{0.5})$ is obtained for HVM and GLS methods for both
Q1 and Q2 elements.
It is consequence of the lack of global regurarity of the exact solution
combined with the use of continuous interpolations to
approximate a discontinuous field.
The global post-processing with interface of discontinuity (GPPID) presents
convergence orders typical of regular solutions: $O(h^{2.0})$ for
Q1 and $O(h^{2.5})$ for Q2. Though not proved, the optimal convergence order
of velocity in $L^2$ norm is usually observed for bilinear elements
\cite{LOULA95}.
The results corresponding to the local post-processing (LPP)
show the optimal orders:
$O(h^{2.0})$ and $O(h^{3.0})$ for Q2, in agreement with the
analysis of homogeneous problems with regular solutions.

\begin{figure}[htb]
\centering
\subfloat[Bilinear.]{\includegraphics[angle=0,width=.49\textwidth,angle=0]{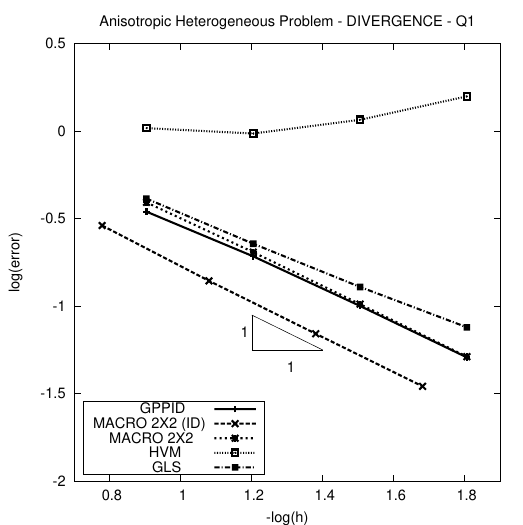}
\label{fig:div-q1}}
%
% \hspace{1cm}
%
\subfloat[Biquadratic.]{\includegraphics[angle=0,width=.49\textwidth,angle=0]{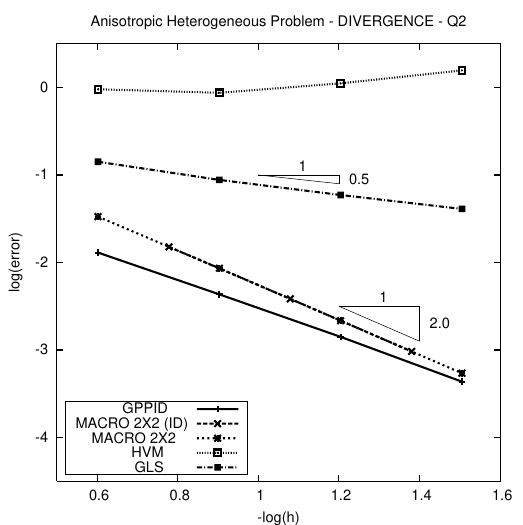}
 \label{fig:div-q2}}
\caption{Convergence of the divergence of the velocity field.}
    \label{fig:bar-conv-divv}%
\end{figure}

Figures \ref{fig:div-q1} and \ref{fig:div-q2} present convergence
results for the divergence of the velocity approximations in
$L^2(\Omega)$-norm.
The poor accuracy of GLS approximation for the velocity on the
interface degradates the orders of convergence of its divergence. The
same orders of convergence obtained for the velocity approximation in
$L^2$-norm, close to $O(h^{0.5})$ for for both Q1 and Q2 elements,
are observed for the divergence in the same norm. This result can be
explained by estimate  (\ref{eq:boundgls}) for the GLS
approximations in ${H(\div)}$-norm for the velocity.
Concerning HVM approximations, according to (\ref{eq:hvm-div}) we
should not expect convergence for divergence using Q1 elements even
for regular exact solutions, except for the ``superconvergence''
usually observed with linear and bilinear $C^0$ velocity
approximations.
But $O(h)$ convergence is expected with Q2 elements for regular
solutions. However, no convergence is observed. This can be
explicated by estimate (\ref{eq:hvm-div}) combined with the lack of
global regularity of the exact solution of this model problem and
the use of  $C^0$ Lagrangian interpolation for the velocity field.

For exactly impose the constraints or naturally accommodate
discontinuities on the interface, the convergence orders observed for
GPPID and LPP are in agreement with those predicted by the numerical
analysis.

\begin{figure}[htb]
\centering
\subfloat[Bilinear.]{\includegraphics[angle=0,width=.49\textwidth,angle=0]{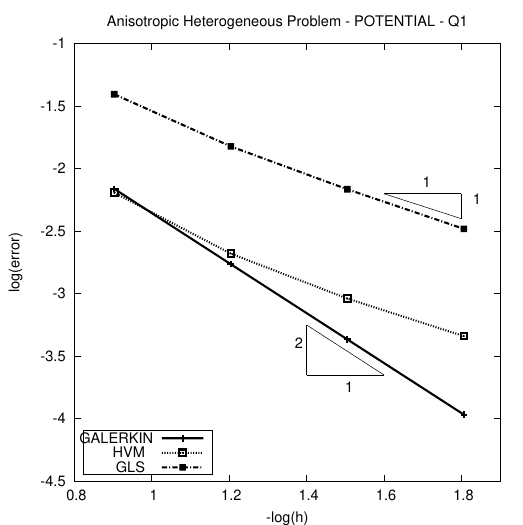}
\label{fig:pot-q1}}
%
% \hspace{1cm}
%
\subfloat[Biquadratic.]{\includegraphics[angle=0,width=.49\textwidth,angle=0]{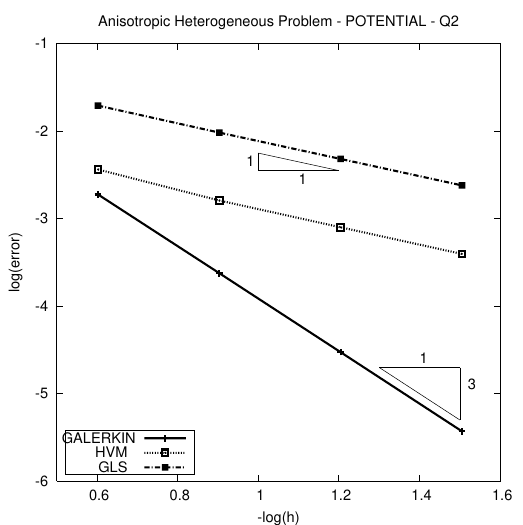}
 \label{fig:pot-q2}}
\caption{Convergence of the potential field.}
    \label{fig:bar-conv-pot}%
\end{figure}

Finally, in Figures \ref{fig:pot-q1} and \ref{fig:pot-q2} we plot
the errors for potential using Q1 and Q2 elements, respectively. The
Galerkin approximation, which is the starting point for the post-processings, shows
the optimal orders predicted by (\ref{eq:errop}), while the
mixed methods show convergence orders of lower order than those we would expected from
(\ref{eq:mpe-p}) for GLS and (\ref{eq:hvm-p}) for HVM.
This can be explained by the coupled estimates for potential and
velocity
 approximations characteristic of mixed methods
((\ref{eq:boundgls}), for GLS, and (\ref{eq:boundhvm}), for HVM) and the lack of global
regularity of the exact solution.

%===================================================================
%
%
\section{CONCLUSIONS}
\label{sec:conclu}
%
%
%===================================================================

Stabilized mixed methods and global post-processings with $C^0$
Lagrangian interpolation present highly stable and accurate
approximations  for continuous velocity fields, as demonstrated and
numerically confirmed in a large number of experiments
\cite{LOULA88,ZIENK92,LOULA95,LOULA99,MASUD2002,LOULA2006}. However,
continuous interpolations are not adequate to approximate
discontinuous velocity fields in heterogeneous porous media with an
interface of material discontinuity.
We proposed generalizations of global and  local post-processing
techniques, based on stabilized variational formulations, for Darcy
flow in heterogeneous porous media with interfaces of
discontinuities. By  exactly imposing the continuity/discontinuity constraints
 on the interface of discontinuity we can  recover  stability,  accuracy and the optimal
 of convergence of classical Lagrangian interpolations.
 Numerical results illustrate the performance of
the proposed methods. Convergence studies for a heterogeneous and
anisotropic porous medium with a smooth interface of discontinuity
confirm the same orders of convergence predicted for homogeneous
problem with smooth solutions, for both global and local
post-processings.
With the local post-processing, the continuity/discontinuity
interface constraints can be exactly imposed in the interior of the
macroelements or naturally accommodated on their edges. In both
situations optimal orders of convergence were obtained.

\section*{Acknowledgements}

The authors thank to  Funda\c c\~ao Carlos Chagas Filho de Amparo \`a
Pesquisa do Estado do Rio de Janeiro (FAPERJ) and to Conselho Nacional
de Desenvolvimento Cient\'{\i}fico e Tecnol\'ogico (CNPq)
%the National Council for Scientific and Technological Development (CNPq)
for the sponsoring.

\bibliographystyle{amsplain}

\end{document}